\documentclass{article}
\usepackage{amssymb}
\usepackage{amsmath}
\usepackage{amsfonts}

\newtheorem{theorem}{Theorem}

\newtheorem{definition}[theorem]{Definition}

\newtheorem{lemma}[theorem]{Lemma}

\newtheorem{remark}[theorem]{Remark}

\numberwithin{equation}{section}
\numberwithin{theorem}{section}

\begin{document}
\title{Regularity in the Obstacle Problem for Parabolic Non-divergence Operators of
H\"{o}rmander type}
\author{Marie Frentz\thanks{email: marie.frentz@gmail.com}\\Department of Mathematics\\London School of Economics\\Houghton Street\\London\\WC2A 2AE\\United Kingdom }
\maketitle
\begin{abstract}
\noindent In this paper we continue the study initiated in \cite{FGN}
concerning the obstacle problem for a class of parabolic nondivergence
operators structured on a set of vector fields $X=\{X_{1},...,X_{q}\}$ in
$\mathbf{R}^{n}$ with $C^{\infty}$-coefficients satisfying H{\"{o}}rmander's
finite rank condition, i.e., the rank of $\mbox{Lie}\lbrack X_{1},...,X_{q}]$
equals $n$ at every point in $\mathbf{R}^{n}$. In \cite{FGN} we proved, under
appropriate assumptions on the operator and the obstacle, the existence and
uniqueness of strong solutions to a general obstacle problem. The main result
of this paper is that we establish further regularity, in the interior as well
as at the initial state, of strong solutions. Compared to \cite{FGN} we in
this paper assume, in addition, that there exists a homogeneous Lie group
$\mathbf{G}=(\mathbf{R}^{n}\!,\circ,\delta_{\lambda})$ such that $X_{1}%
,\ldots,X_{q}$ are left translation invariant on $\mathbf{G}$ and such that
$X_{1},\ldots,X_{q}$ are $\delta_{\lambda}$-homogeneous of degree one.

\vspace{.2in}



\end{abstract}

\section{Introduction}

In this paper we consider the obstacle problem for a class of second order
parabolic subelliptic partial differential equations in non-divergence form
and modeled on a system of vector fields satisfying the H{\"{o}}rmander's
finite rank condition. In particular, we consider operators
\begin{equation}
\mathcal{H}=\sum_{i,j=1}^{q}a_{ij}(x,t)X_{i}X_{j}+\sum_{i=1}^{q}%
b_{i}(x,t)X_{i}-\partial_{t},\label{operator}%
\end{equation}
where $(x,t)\in\mathbf{R}^{n+1}$, $q$ is a positive integer and the functions
$\{a_{ij}(\cdot,\cdot)\}$ and $\{b_{i}(\cdot,\cdot)\}$ are bounded and
measurable on $\mathbf{R}^{n+1}$. In the following, and throughout the paper,
we by $D\subset\mathbf{R}^{n+1}$ denote a bounded and open cylindrical domain
of the form $D=\Omega\times(T_{1},T_{2}]$ where $\Omega\subset\mathbf{R}^{n}$
is a bounded domain and $-\infty<T_{1}<T_{2}<\infty$. We let $\partial_{p}D$
denote the parabolic boundary of $D$. Let $g,f,\varphi:\bar{D}\rightarrow
\mathbf{R}^{n+1}$ be such that $g\geq\varphi$ on $\bar{D}$ and assume that
$g,f,\varphi$ are continuous and bounded on $\bar{D}$. We consider the
following obstacle problem,
\begin{equation}%
\begin{cases}
\max\{\mathcal{H}u(x,t)-f(x,t),\varphi(x,t)-u(x,t)\}=0, & \text{in}\ D,\\
u(x,t)=g(x,t), & \text{on}\ \partial_{p}D.
\end{cases}
\label{e-obs}%
\end{equation}
The purpose of the paper is to advance the mathematical theory for the
obstacle problem in \eqref{e-obs} and in particular to continue the study of
the obstacle problem initiated in \cite{FGN} where a number of important steps
were taken towards developing a rigorous existence theory for the problem in
\eqref{e-obs}. The main result in \cite{FGN} is the existence and uniqueness
of a strong solution to the problem in \eqref{e-obs} in certain bounded
cylindrical domains $D$. Recall that we say that $u\in\mathcal{S}%
_{X,\text{loc}}^{1}(D)\cap C(\overline{D})$ is a strong solution to problem
\eqref{e-obs} if the differential inequality is satisfied a.e. in $D$ and the
boundary data is attained at any point of $\partial_{p}D$. Here $C(\overline
{D})$ is the space of functions, continuous on $\overline{D},$ and the
Sobolev-Stein space $\mathcal{S}_{X,\text{loc}}^{1}(D)$ is defined in Section
2. The subscript $X$ indicates that the space is defined with respect to the
system $X=\{X_{1},...,X_{q}\}$. Moreover, while the study in \cite{FGN}
focused on existence and uniqueness results, the main purpose of this paper is
to establish further regularity for solutions to the problem \eqref{e-obs}. In
this paper we will, compared to \cite{FGN}, impose more restrictive
assumptions on the system $X=\{X_{1},...,X_{q}\}.$ To discuss the structural
assumptions imposed on the operator $\mathcal{H}$, as well as the regularity
assumptions on $a_{ij}$, $b_{i}$, $f$, $\varphi$ and $g$, we first note that
we assume that the system $X=\{X_{1},...,X_{q}\}$ is a set of vector fields in
$\mathbf{R}^{n}$ with $C^{\infty}$-coefficients, i.e., $X_{i}=(c_{i1}%
(x),...,c_{in}(x))\cdot\nabla$, $\nabla=(\partial_{x_{1}},...,\partial_{x_{n}%
})$, where $c_{ik}\in C^{\infty}(\mathbf{R}^{n})$ and $\cdot$ denotes scalar
product in $\mathbf{R}^{n}.$ For technical reasons we will assume that
$c_{ik}$ are bounded; however, this is not restrictive, since the obstacle
problem is stated in a bounded domain. In this paper we impose two essential
restrictions on the system $X=\{X_{1},...,X_{q}\}$. Firstly, we assume that
there exists a homogeneous Lie group $\mathbf{G}=(\mathbf{R}^{n}%
\!,\circ,\delta_{\lambda})$, we refer the reader to Section 2 for the
definition of $\mathbf{G}$, such that
\begin{align}
(i)  &
\mbox{$X_1,\ldots,X_q$ are left translation invariant on $\mathbf{G}$},\nonumber\\
(ii)  &
\mbox{$X_1,\ldots,X_q$ are $\delta_\lambda$-homogeneous of degree one}.\label{1.1x}%
\end{align}
Secondly, we assume that%

\begin{equation}
\text{the vector fields }X=\{X_{1},...,X_{q}\}\text{ satisfy H\"{o}rmanders
finite rank condition.}\label{1.1x+}%
\end{equation}
To be more precise, recall that the Lie-bracket between two vector fields
$X_{i}$ and $X_{j}$ is defined as $[X_{i},X_{j}]=X_{i}X_{j}-X_{j}X_{i}$ and
for an arbitrary multiindex $\theta=(\theta_{1},...,\theta_{d}),\ d\geq1$, we
define%
\[
X^{\theta}=[X_{\theta_{d}},[X_{\theta_{d-1}},...[X_{\theta_{2}},X_{\theta_{1}%
}]]].
\]
The system $X=\{X_{1},...,X_{q}\}$ is said to satisfy a H\"{o}rmander finite
rank condition of order $s$ if $s<\infty$ is the least integer for which the
vector space spanned by $\{X_{i}^{\theta}(x):i=1,...,q,\ |\theta|\leq s\}$ is
$\mathbf{R}^{n},$ for all $x\in\mathbf{R}^{n}$. Note that in \cite{FGN} we did
not assume \eqref{1.1x}. At the end of the paper we discuss to which extent
the assumptions in \eqref{1.1x} are necessary for our main results. Concerning
the $q\times q$ matrix-valued function $A=A(x,t)=\{a_{ij}(x,t)\}=\{a_{ij}\}$
we assume that $A$ is real symmetric with bounded and measurable entries and
that
\begin{equation}
\Lambda^{-1}|\xi|^{2}\leq\sum_{i,j=1}^{q}a_{ij}(x,t)\xi_{i}\xi_{j}\leq
\Lambda|\xi|^{2}\mbox{ whenever }(x,t)\in\mathbf{R}^{n+1},\ \xi\in
\mathbf{R}^{q},\label{1.2}%
\end{equation}
for some $\Lambda$, $1\leq\Lambda<\infty$. Let $d_{p}((x,t),(y,s))=(d(x,y)^{2}%
+|t-s|)^{1/2}$ where $d(x,y)$ is the Carnot-Carath\'{e}odory distance, between
$x,y\in\mathbf{R}^{n}$, induced by $\{X_{1},...,X_{q}\},$ see Section
\ref{secFiol}. Concerning the regularity of $a_{ij}$ and $b_{i}$ we will
assume that $a_{ij}$ and $b_{i}$ have further regularity beyond being only
bounded and measurable. In fact, we assume that
\begin{equation}
a_{ij},\ b_{i}\in\mathcal{C}_{X}^{0,\alpha}(D)\mbox{ whenever }i,j\in
\{1,...,q\},\label{1.2+-}%
\end{equation}
where $\mathcal{C}_{X}^{0,\alpha}(D)$ is the space of functions which are
bounded and H{\"{o}}lder continuous on $D$ and where the H{\"{o}}lder
continuity is defined in terms of the (parabolic) distance induced by the
vector fields. In particular, we assume that there exists a constant $c$,
$0<c<\infty$, such that
\[
|a_{ij}(x,t)-a_{ij}(y,s)|+|b_{i}(x,t)-b_{i}(y,s)|\leq c(d_{p}%
((x,t),(y,s)))^{\alpha},
\]
whenever $(x,t),(y,s)\in D$, $i,j\in\{1,...,q\}$.

We are now ready to formulate the main results proved in this paper. We refer
to Section \ref{SecFS} for the definition of the H\"{o}lder spaces
$\mathcal{C}_{X}^{m,\alpha}$, $m\in\{0,1,2\}$, $\alpha\in(0,1)$, and
Sobolev-Stein spaces $\mathcal{S}_{X}^{p}$, $1\leq p\leq\infty$. When we in
the following write that a constant $c$ depends on the operator $\mathcal{H}$,
$c=c(\mathcal{H})$, we mean that the constant $c$ only depends on $n$, $q$,
$X=\{X_{1},...,X_{q}\}$, $\Lambda$ and $\mathbf{c}_{\alpha}$ where
\[
\mathbf{c}_{\alpha}=\sum_{i,j=1}^{q}\Vert a_{ij}\Vert_{\mathcal{C}%
_{X}^{0,\alpha}(D)}+\sum_{i=1}^{q}\Vert b_{i}\Vert_{\mathcal{C}_{X}^{0,\alpha
}(D)}.
\]
Given bounded domains $\Omega,\Omega^{\prime}$ in $\mathbf{R}^{n}$, and $T>0$,
we let $\Omega_{T}=\Omega\times(0,T]$, $\Omega_{T}^{\prime}=\Omega^{\prime
}\times(0,T]$. We prove the following two theorems.

\begin{theorem}
\label{th-int} Let $\mathcal{H}$ be defined as in \eqref{operator}, assume
\eqref{1.1x}-\eqref{1.2+-}, let $\Omega,\ \Omega^{\prime}$ be bounded domains
in $\mathbf{R}^{n}$ such that $\Omega^{\prime}\subset\subset\Omega$ and let
$0<T^{\prime\prime}<T^{\prime}<T$. Let $g,f,\varphi:\bar{\Omega}%
_{T}\rightarrow\mathbf{R}^{n+1}$ be such that $g\geq\varphi$ on $\bar{\Omega
}_{T}$ and assume that $g,f,\varphi$ are continuous and bounded on
$\bar{\Omega}_{T}$. Let $\alpha\in(0,1)$ and let $u$ be a strong solution to
problem \eqref{e-obs} in $\Omega_{T}$. Then the following holds:

\begin{itemize}
\item[i)] if $\varphi\in\mathcal{C}_{X}^{0,\alpha}(\Omega_{T})$ then
$u\in\mathcal{C}_{X}^{0,\alpha}(\Omega^{\prime}\times(T^{\prime\prime
},T^{\prime}))$ and
\[
\Vert u\Vert_{\mathcal{C}_{X}^{0,\alpha}(\Omega^{\prime}\times(T^{\prime
\prime},T^{\prime}))}\leq c\left(  \alpha,\Omega,\Omega^{\prime},T,T^{\prime
},T^{\prime\prime},\mathcal{H},\Vert f\Vert_{\mathcal{C}_{X}^{0,\alpha}%
(\Omega_{T})},\Vert g\Vert_{L^{\infty}(\Omega_{T})},\Vert\varphi
\Vert_{\mathcal{C}_{X}^{0,\alpha}(\Omega_{T})}\right)  ;
\]

\item[ii)] if $\varphi\in\mathcal{C}_{X}^{1,\alpha}(\Omega_{T})$ then
$u\in\mathcal{C}_{X}^{1,\alpha}(\Omega^{\prime}\times(T^{\prime\prime
},T^{\prime}))$ and
\[
\Vert u\Vert_{\mathcal{C}_{X}^{1,\alpha}(\Omega^{\prime}\times(T^{\prime
\prime},T^{\prime}))}\leq c\left(  \alpha,\Omega,\Omega^{\prime},T,T^{\prime
},T^{\prime\prime},\mathcal{H},\Vert f\Vert_{\mathcal{C}_{X}^{0,\alpha}%
(\Omega_{T})},\Vert g\Vert_{L^{\infty}(\Omega_{T})},\Vert\varphi
\Vert_{\mathcal{C}_{X}^{1,\alpha}(\Omega_{T})}\right)  ;
\]

\item[iii)] if $\varphi\in\mathcal{C}_{X}^{2,\alpha}(\Omega_{T})$ then
$u\in\mathcal{S}_{X}^{\infty}(\Omega^{\prime}\times(T^{\prime\prime}%
,T^{\prime}))$ and
\[
\Vert u\Vert_{\mathcal{S}_{X}^{\infty}(\Omega^{\prime}\times(T^{\prime\prime
},T^{\prime}))}\leq c\left(  \alpha,\Omega,\Omega^{\prime},T,T^{\prime
},T^{\prime\prime},\mathcal{H},\Vert f\Vert_{\mathcal{C}_{X}^{0,\alpha}%
(\Omega_{T})},\Vert g\Vert_{L^{\infty}(\Omega_{T})},\Vert\varphi
\Vert_{\mathcal{C}_{X}^{2,\alpha}(\Omega_{T})}\right)  .
\]

\end{itemize}
\end{theorem}

\begin{theorem}
\label{th} Let $\mathcal{H}$ be defined as in \eqref{operator}, assume
\eqref{1.1x}-\eqref{1.2+-}, let $\Omega,\ \Omega^{\prime}$ be bounded domains
in $\mathbf{R}^{n}$ such that $\Omega^{\prime}\subset\subset\Omega$ and let
$0<T^{\prime}<T$. Let $g,f,\varphi:\bar{\Omega}_{T}\rightarrow\mathbf{R}%
^{n+1}$ be such that $g\geq\varphi$ on $\bar{\Omega}_{T}$ and assume that
$g,f,\varphi$ are continuous and bounded on $\bar{\Omega}_{T}$. Let $\alpha
\in(0,1)$ and let $u$ be a strong solution to problem \eqref{e-obs} in
$\Omega_{T}$. Then the following holds:

\begin{itemize}
\item[i)] if $g,\varphi\in\mathcal{C}_{X}^{0,\alpha}(\Omega_{T})$ then
$u\in\mathcal{C}_{X}^{0,\alpha}(\Omega_{T^{\prime}}^{\prime})$ and
\[
\Vert u\Vert_{\mathcal{C}_{X}^{0,\alpha}(\Omega_{T^{\prime}}^{\prime})}\leq
c\left(  \alpha,\Omega_{T},\Omega_{T\prime}^{\prime},\mathcal{H},\Vert
f\Vert_{\mathcal{C}_{X}^{0,\alpha}(\Omega_{T})},\Vert g\Vert_{\mathcal{C}%
_{X}^{0,\alpha}(\Omega_{T})},\Vert\varphi\Vert_{\mathcal{C}_{X}^{0,\alpha
}(\Omega_{T})}\right)  ;
\]

\item[ii)] if $g,\varphi\in\mathcal{C}_{X}^{1,\alpha}(\Omega_{T})$ then
$u\in\mathcal{C}_{X}^{1,\alpha}(\Omega_{T^{\prime}}^{\prime})$ and
\[
\Vert u\Vert_{\mathcal{C}_{X}^{1,\alpha}(\Omega_{T^{\prime}}^{\prime})}\leq
c\left(  \alpha,\Omega_{T},\Omega_{T\prime}^{\prime},\mathcal{H},\Vert
f\Vert_{\mathcal{C}_{X}^{0,\alpha}(\Omega_{T})},\Vert g\Vert_{\mathcal{C}%
_{X}^{1,\alpha}(\Omega_{T})},\Vert\varphi\Vert_{\mathcal{C}_{X}^{1,\alpha
}(\Omega_{T})}\right)  ;
\]

\item[iii)] if $g,\varphi\in\mathcal{C}_{X}^{2,\alpha}(\Omega_{T})$ then
$u\in\mathcal{S}_{X}^{\infty}(\Omega_{T^{\prime}}^{\prime})$ and
\[
\Vert u\Vert_{\mathcal{S}_{X}^{\infty}(\Omega_{T^{\prime}}^{\prime})}\leq
c\left(  \alpha,\Omega_{T},\Omega_{T\prime}^{\prime},\mathcal{H},\Vert
f\Vert_{\mathcal{C}_{X}^{0,\alpha}(\Omega_{T})},\Vert g\Vert_{\mathcal{C}%
_{X}^{2,\alpha}(\Omega_{T})},\Vert\varphi\Vert_{\mathcal{C}_{X}^{2,\alpha
}(\Omega_{T})}\right)  .
\]

\end{itemize}
\end{theorem}

Theorem \ref{th-int} concerns the optimal interior regularity for the solution
$u$ to the obstacle problem under different assumptions on the regularity of
the obstacle $\varphi$. In particular, Theorem \ref{th-int} treats the case of
non-smooth obstacles as well as the case of smooth obstacles. The results
stated in the theorems are similar: the solution is, up to $\mathcal{S}%
_{X}^{\infty}$-smoothness, as smooth as the obstacle. Note that Theorem
\ref{th} gives similar results but in this case $\Omega_{T^{\prime}}^{\prime}$
is not a compact subset of $\Omega_{T}$. As such Theorem \ref{th} concerns the
optimal regularity up to the initial state for the solution $u$ to the
obstacle problem under the stated assumptions on the regularity of the
obstacle $\varphi$. Based on the discussion below we claim that Theorem
\ref{th-int} and Theorem \ref{th} represent two new contributions to the
literature on optimal regularity for obstacle problems. In fact, we are not
aware of any previous works, in the parabolic and genuinely subelliptic case,
i.e., when $X=\{X_{1},...,X_{q}\}$ is not identical to $\{\partial_{x_{1}%
},...,\partial_{x_{n}}\}$, devoted to the obstacle problem in \eqref{e-obs}.

To outline previous work on optimal interior regularity in obstacle problems
we first discuss the more classical case when $X=\{X_{1},...,X_{q}%
\}\equiv\{\partial_{x_{1}},...,\partial_{x_{n}}\}$ and we note that in this
case there is an extensive literature on the existence of generalized
solutions to the obstacle problem in Sobolev spaces starting with the
pioneering papers \cite{McK}, \cite{vM}, \cite{vM1} and \cite{F75}.
Furthermore, optimal regularity of the solution to the obstacle problem for
the Laplace equation was first proved by Caffarelli and Kinderlehrer
\cite{CK80} and we note that the techniques used in \cite{CK80} are based on
the Harnack inequality for harmonic functions and the control of a harmonic
function by its Taylor expansion. The most extensive and complete treatment of
the obstacle problem for the heat equation can be found in Caffarelli,
Petrosyan and Shahgholian \cite{CPS} and it is interesting to note that most
of the arguments in \cite{CPS} make use of a blow-up technique previously also
used by Caffarelli, Karp and Shahgholian in \cite{CKS00} in the stationary
case. We here also mention the paper \cite{BDM06} where the optimal regularity
of the obstacle problem for second order uniformly elliptic parabolic
equations has been proved by a method inspired by the original one in
\cite{CK80} based on the Harnack inequality. On the other hand the blow-up
method has been employed in more general settings in \cite{PS07},
\cite{Sha08}. Concerning previous work on the optimal interior regularity in
obstacle problems in a subelliptic setting we note that in \cite{DGS03} the
obstacle problem is considered for the strongly degenerate case of
sublaplacians on Carnot groups. The paper \cite{DGP07} addresses, in the same
framework, the study of the regularity of the free boundary. In particular,
the sublaplacian considered in \cite{DGS03} can be considered as a special
case of the stationary versions of the more general operators studied in this
paper. Finally, we note that in \cite{FNPP} the author, together with
Nystr\"{o}m, Pascucci and Polidoro, recently have established an appropriate
version of Theorem \ref{th-int} for a class of second order differential
operators of Kolmogorov type of the form
\begin{equation}
H=\sum_{i,j=1}^{m}a_{ij}(x,t)\partial_{x_{i}x_{j}}+\sum_{i=1}^{m}%
b_{i}(x,t)\partial_{x_{i}}+\sum_{i,j=1}^{n}b_{ij}x_{i}\partial_{x_{j}%
}-\partial_{t}\label{1.1ahhha}%
\end{equation}
where $(x,t)\in\mathbf{R}^{n+1}$, $m$ is a positive integer satisfying $m\leq
n$, the functions $\{a_{ij}(\cdot,\cdot)\}$ and $\{b_{i}(\cdot,\cdot)\}$ are
continuous and bounded and the matrix $B=\{b_{ij}\}$ is a matrix of constant
real numbers. The structural assumptions imposed in \cite{FNPP} on the
operator $H$ implies that $H$ is a hypoelliptic ultraparabolic operator of
Kolmogorov type. Note however, that the operator in \eqref{1.1ahhha} is
different from the class of operators considered in this paper due to the fact
that space and time couple through the lower order term $Y=\sum_{i,j=1}%
^{n}b_{ij}x_{i}\partial_{x_{j}}-\partial_{t}$ and since $\{X_{1}%
,...,X_{m}\}=\{\partial_{x_{1}},...,\partial_{x_{m}}\}$.

To outline previous work on the optimal regularity up to the initial state in
obstacle problems we note that there is, already in the case when
$X=\{X_{1},...,X_{q}\}\equiv\{\partial_{x_{1}},...,\partial_{x_{n}}\}$, a very
limited literature on this topic. In fact, in this case we are only aware of
the results by Nystr{\"{o}}m \cite{N}, Shahgholian \cite{Sha08} (see also
Petrosyan and Shahgholian \cite{PS07}). While the arguments in \cite{Sha08}
allow for certain fully non-linear parabolic equations, in \cite{N} the
techniques were conveyed in the context of pricing of multi-dimensional
American options in a financial market driven by a general multi-dimensional
It\^{o} diffusion. In \cite{N} the machinery and techniques were developed and
described assuming more regularity on the operator and the obstacle than
needed and in the standard context of American options. However, the results
in \cite{N} and \cite{Sha08} do not apply, for example and considering
financial applications, in the setting of Asian options or the Hobson-Rogers
model for stochastic volatility \cite{HR98}. Therefore\ Nystr\"{o}m, Pascucci
and Polidoro, in \cite{NPP}, recently established an appropriate version of
Theorem \ref{th} for the class of second order differential operators of
Kolmogorov type briefly defined in \eqref{1.1ahhha}, and we note that
\cite{NPP} is a continuation of the study initialized in \cite{FNPP}. We also
note that the results in \cite{NPP} also apply to uniformly parabolic
equations, i.e., the case when $X=\{X_{1},...,X_{q}\}\equiv\{\partial_{x_{1}%
},...,\partial_{x_{n}}\}$. In this case the results in \cite{NPP} slightly
improve upon Theorem 4.3 in \cite{PS07} (see also Theorems 1.2 and 1.3 in
\cite{Sha08}) since in \cite{NPP} the authors get the H\"{o}lder regularity of
the solution with the optimal exponent.

Concerning our proofs of Theorem \ref{th-int} and Theorem \ref{th} we will
proceed in a fashion structurally similar to the corresponding proofs in
\cite{FNPP} and \cite{NPP}. In particular, we will use the type of blow-up
technique introduced by Caffarelli, Karp and Shahgholian in \cite{CKS00} in
the stationary case and by Caffarelli, Petrosyan and Shahgholian \cite{CPS} in
the study of the heat equation. This method is flexible enough to give
stream-lined proofs of the statements in Theorem \ref{th-int} and Theorem
\ref{th} for non-smooth as well as smooth obstacles. To briefly outline the
method we let $d(x,y)$ be the Carnot-Carath\'{e}odory distance between
$x,y\in\mathbf{R}^{n}$, induced by $\{X_{1},...,X_{q}\}$ and we let
$B_{d}(x,r)=\{y\in\mathbf{R}^{n}:\ d(x,y)<r\}$, whenever $x\in\mathbf{R}^{n} $
and $r>0$, denote the\ balls associated to $d$. We let
\begin{align}
C_{r}(x,t)  & =B_{d}(x,r)\times(t-r^{2},t+r^{2}),\nonumber\\
C_{r}^{+}(x,t)  & =B_{d}(x,r)\times(t,t+r^{2}),\nonumber\\
C_{r}^{-}(x,t)  & =B_{d}(x,r)\times(t-r^{2},t),\label{pol3aa}%
\end{align}
whenever $(x,t)\in\mathbf{R}^{n+1}$ and $r>0$. In fact, we will use a
modification of these cylinders, see the discussion below Lemma \ref{weakmax}.
Using the notation in \eqref{pol3aa} we in this paper build the core part of
the argument, in the case of Theorem \ref{th-int}, on the function
\begin{equation}
S_{k}^{-}(u)=\sup_{C_{2^{-k}}^{-}(0,0)}|u|\label{1.17x}%
\end{equation}
and, in the case of Theorem \ref{th}, on the function
\begin{equation}
S_{k}^{+}(u)=\sup_{C_{2^{-k}}^{+}(0,0)}|u|.\label{1.17x+}%
\end{equation}
In \eqref{1.17x} and \eqref{1.17x+} $u$ is a solution to the obstacle problem
in $C_{1}^{-}(0,0)$ and $C_{1}^{+}(0,0)$, respectively. We will assume that
the quadruple $(u,f,g,\varphi)$, which specifies the degrees of freedom in the
obstacle problem, belongs to certain function classes defined in Subsection
\ref{SecFC}. Moreover given $\varphi$, in this construction we let, in the
proof of Theorem \ref{th-int}, $F$ and $\gamma$ be determined as follows:
\begin{align}
(i)  & F=P_{0}^{(0,0)}\varphi,\ \gamma=\alpha,\nonumber\\
(ii)  & F=P_{1}^{(0,0)}\varphi,\ \gamma=1+\alpha,\nonumber\\
(iii)  & F=\varphi,\ \gamma=2,\label{lex1}%
\end{align}
where $\alpha\in(0,1)$ and $P_{m}^{(0,0)}\varphi$ is a certain intrinsic
Taylor expansion associated to $\varphi$ and as outlined in Subsection
\ref{SecFS}. In particular, as an important step in the proof of Theorem
\ref{th-int}, we prove that there exists a positive constant $c$ such that,
for all $k\in\mathbf{N}$,
\begin{equation}
S_{k+1}^{-}(u-F)\leq\max\biggl (c\,2^{-(k+1)\gamma},\frac{S_{k}^{-}%
(u-F)}{2^{\gamma}},\frac{S_{k-1}^{-}(u-F)}{2^{2\gamma}},\dots,\frac{S_{0}%
^{-}(u-F)}{2^{(k+1)\gamma}}\biggr ).\label{1.17z}%
\end{equation}
In the case of Theorem \ref{th}, \ we again assume that $(u,f,g,\varphi)$
belongs to certain function classes and given $g$, in this construction, we
let $F$ and $\gamma$ be determined as follows:
\begin{align}
(i)  & F=P_{0}^{(0,0)}g,\ \gamma=\alpha,\nonumber\\
(ii)  & F=P_{1}^{(0,0)}g,\ \gamma=1+\alpha,\nonumber\\
(iii)  & F=g,\ \gamma=2.\label{lex2}%
\end{align}
Similarly, to prove Theorem \ref{th}, we prove that there exists a positive
constant $c$ such that, for all $k\in\mathbf{N}$,
\begin{equation}
S_{k+1}^{+}(u-F)\leq\max\biggl (c\,2^{-(k+1)\gamma},\frac{S_{k}^{+}%
(u-F)}{2^{\gamma}},\frac{S_{k-1}^{+}(u-F)}{2^{2\gamma}},\dots,\frac{S_{0}%
^{+}(u-F)}{2^{(k+1)\gamma}}\biggr ).\label{1.17z+}%
\end{equation}
Note that we are able to carry the argument through using only $g$, instead of
$\varphi$, in \eqref{lex2} and this is one key aspect of the argument. In
either case the proofs of \eqref{1.17z} and \eqref{1.17z+} are based on
arguments by contradiction and the argument in the case of \eqref{1.17z+}
differs at key points compared to the corresponding proof in the case
\eqref{1.17z} due to the presence, in the case \eqref{1.17z+}, of the boundary
at $t=0$.

The rest of the paper is organized as follows. Section 2 is of preliminary
nature and we here introduce relevant notations and recall some basic facts
concerning homogeneous Lie groups. We also introduce function spaces and
consider their characterizations based on polynomials. In Section 3 we state
some results concerning operators $\mathcal{H}$ defined as in \eqref{operator}
assuming \eqref{1.1x}-\eqref{1.2+-}. These results concern existence and
estimates for fundamental solutions, the Harnack inequality and Schauder
estimates. In this section we also use these results to prove a few auxiliary
technical estimates to be used in the proofs of Theorem \ref{th-int} and
Theorem \ref{th}. In Section 4 we derive certain estimates for the obstacle
problem to be used in the proofs of Theorem \ref{th-int} and Theorem \ref{th}.
Section 5 is devoted to the proof of Theorem \ref{th-int} while Section 6 is
devoted to the proof of Theorem \ref{th}. Finally, in Section 7 we complete
this paper with some concluding remarks.

\section{Preliminaries}

In this section we firstly introduce some notations and recall some basic
notions concerning homogeneous Lie groups. We refer to the monograph
\cite{BLU07} for a detailed treatment of the subject. Secondly we introduce
relevant function spaces and consider their Taylor approximation. Thirdly, and
finally, we introduce certain function classes related to the Cauchy-Dirichlet
problem and to the obstacle problem. Based on these function classes we are
able to give streamlined proofs of important lemmas and to prove Theorem
\ref{th-int} and Theorem \ref{th} in a unified way.

\subsection{Homogeneous Lie groups\label{secFiol}}

Let $\circ$ be a given group law on $\mathbf{R}^{n}$ and suppose that the map
$(x,y)\mapsto y^{-1}\circ x$ is smooth. Then $\mathbf{G}=(\mathbf{R}^{n}%
,\circ)$ is called a \emph{Lie group}. On $\mathbf{G}$ we define the left
translation operator $\tau_{\alpha}(x)=\alpha\circ x$ and we say that a vector
field $X$ on $\mathbf{R}^{n}$ is \emph{left-invariant} if $X(\tau_{\alpha
}(x))=\tau_{\alpha}(X(x))$. We let $\mathbf{g}$ denote the Lie algebra of
$\mathbf{G}$, i.e., the Lie-algebra of left-invariant vector fields on
$\mathbf{R}^{n}$. Moreover, $\mathbf{G}$ is said to be \emph{homogeneous} if
there exists a family of \emph{dilations} $\left(  \delta_{\lambda}\right)
_{\lambda>0}$ on $\mathbf{G}$ of the form
\begin{equation}
\delta_{\lambda}(x)=\delta_{\lambda}(x^{(1)},...,x^{(l)})=(\lambda
x^{(1)},...,\lambda^{l}x^{(l)})=(\lambda^{\sigma_{1}}x_{1},...,\lambda
^{\sigma_{n}}x_{n}),\label{str1}%
\end{equation}
where $1\leq\sigma_{1}\leq...\leq\sigma_{n},$ and if $\left(  \delta_{\lambda
}\right)  _{\lambda>0}$ defines an automorphism of the group, \emph{i.e.},
\[
\delta_{\lambda}(x\circ y)=\left(  \delta_{\lambda}(x)\right)  \circ\left(
\delta_{\lambda}(y)\right)  ,\quad\text{for all}\ x,y\in\mathbf{R}%
^{n}\ \text{and}\ \lambda>0.
\]
Note that in \eqref{str1} we have that $x^{(i)}\in\mathbf{R}^{n_{i}}$ for
$i\in\{1,...,l\}$ and $n_{1}+....+n_{l}=n$. In particular, the dilation
$\delta_{\lambda}$ induces a direct sum decomposition on $\mathbf{R}^{n}$
\[
\mathbf{R}^{n}=V_{1}\oplus\dots\oplus V_{l}.
\]
The natural number
\begin{equation}
Q:=\text{dim}V_{1}+2\,\text{dim}V_{2}+\dots+l\text{dim}V_{l}=n_{1}%
+2n_{2}+...+ln_{l},\label{e-Q}%
\end{equation}
is called the \emph{homogeneous dimension} of $\mathbf{G}$ with respect to
$\delta_{\lambda}$. For short we write $(\mathbf{R}^{n},\circ,\delta_{\lambda
})$ for the Lie group $\mathbf{G}$ equipped with the family of dilations
$\left(  \delta_{\lambda}\right)  _{\lambda>0}$. We next recall some useful
facts about homogeneous Lie groups and we first note that \eqref{1.1x} and
\eqref{1.1x+} imply that $n_{1}=q$ and that span$\big\{X_{1}(0),\dots
,X_{q}(0)\big\}=V_{1}$. Hence it is not restrictive to assume $q=\text{dim}%
V_{1}$ and $X_{i}(0)=\mathbf{e}_{i}$ for $i=1,\dots,q$, where $\{\mathbf{e}%
_{i}\}_{1\leq i\leq n}$ denotes the canonical basis of $\mathbf{R}^{n}$. In
particular, $X_{i}$ is, for $i=1,\dots,q$, the unique vector field in
$\mathbf{g}$ that agrees with $\partial/\partial x_{i}$ at the origin. With
these hypotheses $\mathbf{G}=(\mathbf{R}^{n},\circ,\delta_{\lambda})$ can be
referred to as a \emph{homogeneous Carnot group}. We say that $\mathbf{G}$ is
of step $l$ and that it has $q=n_{1}$ generators. In the literature, see,
e.g., \cite{Fo}, \cite{HK}, \cite{RS}, \cite{VSC}, a Carnot group (or a
stratified group) $\mathbf{G}$ is defined as a connected and simply connected
Lie group whose Lie algebra $\mathbf{\tilde{g}}$ admits a stratification
\[
\mathbf{\tilde{g}}=W_{1}\oplus\dots\oplus W_{l}%
\mbox{ with $[W_1,W_i]=W_{i+1}$ and
$[W_1,W_l]=\{0\}$}.
\]
In fact, any homogeneous Carnot group is a Carnot group according to the
classical definition. On the other hand, up to isomorphism, the opposite
implication is also true, see, e.g., \cite{BU04}. To continue we set
\[
||x||_{{\mathbf{G}}}=\left(  \sum_{j=1}^{l}\sum_{i=1}^{m_{j}}\Big(x_{i}%
^{(j)}\Big)^{\frac{2l!}{j}}\right)  ^{\frac{1}{2l!}},
\]
and we observe that the above function is homogeneous of degree $1$ on
$\mathbf{R}^{n}$ in the sense that
\[
\left\Vert \left(  \lambda x^{(1)},....,\lambda^{l}x^{(l)}\right)  \right\Vert
_{{\mathbf{G}}}=\lambda||x||_{{\mathbf{G}}},
\]
for every $x\in\mathbf{R}^{n}$ and for any $\lambda>0$. In particular,
$||\cdot||_{\mathbf{G}}$ is a homogeneous norm on $\mathbf{G,}$ smooth away
from the origin. Using this homogeneous norm we define a quasi-distance
$\hat{d}:\mathbf{R}^{n}\times\mathbf{R}^{n}\rightarrow\mathbf{[}%
0\mathbf{,}\infty)$ by setting $\hat{d}(x,y):=||y^{-1}\circ x||_{\mathbf{G}}$.
The term quasi-distance is used to indicate that $\hat{d}$ has the same
properties as a distance, except that $\hat{d}$ satisfy the quasi-triangle
inequality and not the triangle inequality. We also recall that for every
compact set $K\subset\mathbf{R}^{n}$ there exist two positive constants
$c_{K}^{-}$ and $c_{K}^{+}$, such that
\begin{equation}
c_{K}^{-}|x-y|\leq\hat{d}(x,y)\leq c_{K}^{+}|x-y|^{\frac{1}{s}},\qquad
\text{for all}\ x,y\in K.\label{e-dist-eucl}%
\end{equation}
Here $|\cdot|$ denotes the Euclidean distance and $s$ is the order in
H\"{o}rmander's finite rank condition for $X=\{X_{1},...,X_{q}\}$. For a proof
of \eqref{e-dist-eucl} we refer to \cite{NSW85}. Finally, we recall the notion
of \emph{Carnot-Carath\'{e}odory distance}. In particular, let $X=\{X_{1}%
,...,X_{q}\}$ be a system satisfying \eqref{1.1x+} and let $T>0$. We say that
an absolutely continuous curve $\gamma:[0,T]\rightarrow\mathbf{R}^{n}$ is a
sub-unit curve with respect to the system $X=\{X_{1},...,X_{q}\} $, $\gamma$
is $X$-subunit for short, if there exist measurable functions $h=(h_{1}%
,...,h_{q})$ such that
\[
\gamma^{\prime}(t)=\sum_{j=1}^{q}h_{j}(t)X_{j}(\gamma(t))\ \ \text{for a.e.
$t\in\lbrack0,T]$ with \ }\sum_{j=1}^{q}h_{j}(t)^{2}\leq1\ \ \text{a.e.}%
\]
In the following we let $l(\gamma)$ denote this $T$ and for $x,y\in
\mathbf{R}^{n}$ we define%
\[
d(x,y)=\inf\left\{  l(\gamma)|\gamma:[0,l(\gamma)]\rightarrow\mathbf{R}%
^{n},\ \mbox{$\gamma$ is $X$\text{-subunit}},\ \gamma(0)=x\text{ and }%
\gamma(l(\gamma))=y\right\}  .
\]
It is well known by Chow-Rashevsky's connectivity theorem, see \cite{C} and
\cite{Ras}, that if $X=\{X_{1},...,X_{q}\}$ is a system satisfying
\eqref{1.1x+}, then the above set is nonempty and hence $d(x,y)$ is finite for
every pair of points. In particular, $d$ is referred to as the
Carnot-Carath\'{e}odory distance, $CC$-distance for short, between the points
$x$ and $y$. Moreover, $d$ is a metric on $\mathbf{R}$ and $||x||:=d(0,x)$ is
a symmetric, homogeneous norm on $\mathbf{G}$. Comparing $\hat{d}$ and $d$ one
can prove that there exist positive constants $c_{1},\ c_{2}$ such that
\begin{equation}
c_{1}\hat{d}(x,y)\leq d(x,y)\leq c_{2}\hat{d}(x,y)\mbox{ for all $x,y\in
\mathbf R^n.$}\label{ekvidist}%
\end{equation}
In fact, both $\hat{d}$ and $d$ are homogeneous of degree $1$ with respect to
$\delta_{\lambda}$ and on a homogeneous group all homogeneous norms are
equivalent, see for instance Chapter 5 in \cite{BLU07}. We next develop the
necessary parabolic notation. In particular, given $\mathbf{G}=(\mathbf{R}%
^{n},\circ,\delta_{\lambda})$ as above we extend this Lie group to a Lie group
on $\mathbf{R}^{n+1}$ by defining
\[
(x,t)^{-1}\circ(y,s)=(x^{-1}\circ y,s-t)
\]
whenever $(x,t),(y,s)\in\mathbf{R}^{n+1}$. Furthermore, we extend
$\delta_{\lambda}$ to $\mathbf{R}^{n+1}$ by defining $\delta_{\lambda
}(x,t)=(\delta_{\lambda}(x),\lambda^{2}t)$ and we write $\mathbf{L}$ for the
constructed homogeneous Lie group $(\mathbf{R}^{n+1},\circ,\delta_{\lambda})$.
Notationally we do not differentiate between dilations on $\mathbf{G}$ and
$\mathbf{L}$ respectively, however, it will be clear from the context which
type of dilations that is used. Moreover, we set
\[
||(x,t)||_{p}=d_{p}((0,0),(x,t)),\text{ \ \ }\Vert(x,t)\Vert_{{\mathbf{L}}%
}=\left(  ||x||_{{\mathbf{G}}}^{2\cdot l!}+|t|^{l!}\right)  ^{\frac{1}{2\cdot
l!}},
\]
and we observe that these functions are homogeneous of degree $1$ on
$\mathbf{R}^{n+1}$, in the sense that
\[
\Vert\delta_{\lambda}(x,t)\Vert_{p}=\lambda\Vert(x,t)\Vert_{p},\text{
\ \ }\Vert\delta_{\lambda}(x,t)\Vert_{{\mathbf{L}}}=\lambda\Vert
(x,t)\Vert_{{\mathbf{L}}},
\]
for every $(x,t)\in\mathbf{R}^{n+1}$ and for any $\lambda>0$. Finally, we
extend $\hat{d}$ and $d$ to $\mathbf{R}^{n+1}\times\mathbf{R}^{n+1}$ by
defining
\begin{equation}
d_{p}((x,t),(y,s))=\sqrt{d(x,y)^{2}+|t-s|},\ \ \ \ \hat{d}_{p}%
((x,t),(y,s))=\sqrt{\hat{d}(x,y)^{2}+|t-s|}\ \label{dp}%
\end{equation}
for $(x,t),\ (y,s)\in\mathbf{R}^{n+1}.$ We continue and introduce the
appropriate balls. For $d$ we define $B_{d}(x,r):=\{y\in\mathbf{R}%
^{n}:\ d(x,y)<r\}$, whenever $x\in\mathbf{R}^{n}$ and $r>0$ and similarly, for
$d_{p}$, we define $B_{d_{p}}((x,t),r):=\{(y,s)\in\mathbf{R}^{n+1}%
:d_{p}((x,t),(y,s))<r\}$. Based on $B_{d}(x,r)$ we introduce the bounded
cylindrical domains $C_{r}(x,t)$, $C_{r}^{+}(x,t)$, $C_{r}^{-}(x,t)$ as in
\eqref{pol3aa}. Then, using the homogeneity of the group $\mathbf{G}$ we see
that there exists a constant $c$ such that
\begin{equation}
c^{-1}r^{Q}\leq|B_{d}(x,r)|\leq cr^{Q}\text{ \ \ and \ \ }\ c^{-1}r^{Q+2}%
\leq|C_{r}(x,t)|\leq cr^{Q+2}\label{bvolym}%
\end{equation}
whenever $(x,t)\in\mathbf{R}^{n+1}$ and $r>0$. Here $|B_{d}(x,r)|$ and
$|C_{r}(x,t)|$ denote the Euclidean volume of $B_{d}(x,r)$ and $C_{r}(x,t)$
respectively and $Q$ is the homogeneous dimension of $\mathbf{G}$ as
introduced in \eqref{e-Q}. Throughout the paper we will often write $C$,
$C^{+}$ and $C^{-}$ for the cylinders $C_{1}(0,0)$, $C_{1}^{+}(0,0)$,
$C_{1}^{-}(0,0)$ introduced in \eqref{pol3aa}. Similarly, we will often write
$C_{r}$, $C_{r}^{+}$ and $C_{r}^{-}$ for the cylinders $C_{r}(0,0)$,
$C_{r}^{+}(0,0)$, $C_{r}^{-}(0,0)$.

Finally we define the exponential map $Exp:\mathbf{g}\rightarrow\mathbf{G}$.
For smooth vector fields $X,$ there exists a unique solution to
\[
\left\{
\begin{array}
[c]{l}%
\gamma^{\prime}(t)=XI(\gamma(t))\\
\gamma(0)=x
\end{array}
\right.
\]
where $I$ is the identity map. We denote this solution $\gamma_{X}(t)$ and we
define $\exp(t\cdot X)(x):=\gamma_{X}(t).$ Then $\exp$ is well defined for all
$X\in$ $\mathbf{g}$, $x\in\mathbf{R}^{n}$ and $t\in\mathbf{R}$. The
exponential map $Exp:$ $\mathbf{g}\rightarrow\mathbf{G}$ is defined through%
\[
Exp(X):=\exp(1\cdot X)(0).
\]

\subsection{Function spaces and polynomial approximation\label{SecFS}}

Let $D\subset\mathbf{R}^{n+1}$ be an open domain. We say that $D$ is a bounded
cylindrical domain if
\begin{equation}
\mbox{$D=\Omega\times (T_1,T_2]$ for some $\Omega\subset\mathbf R^n$, open
and bounded, and some $-\infty<T_1<T_2<\infty$}.\label{cyl}%
\end{equation}
Let $\alpha\in(0,1)$ and let $D\subset\mathbf{R}^{n+1}$ be a bounded
cylindrical domain in the sense of (\ref{cyl}). Let $\mathbf{G}=(\mathbf{R}%
^{n},\circ,\delta_{\lambda})$ and $X=\{X_{1},...,X_{q}\}$ be as in
\eqref{1.1x} and \eqref{1.1x+}. Let $\mathbf{L}$ be the parabolic extension of
$\mathbf{G}$ as defined in the previous subsection and let $d_{p}$ be the
parabolic CC-distance introduced in \eqref{dp}. In the following we will often
denote points in $\mathbf{R}^{n+1}$ by $z$ and $\zeta$, i.e., $z=(x,t)$,
$\zeta=(y,s)$. Using this notation we in the following by $\mathcal{C}%
_{X}^{0,\alpha}(D)$, $\mathcal{C}_{X}^{1,\alpha}(D)$ and $\mathcal{C}%
_{X}^{2,\alpha}(D)$ denote H\"{o}lder spaces defined by the following norms:

\begin{eqnarray*}
\Vert u\Vert_{\mathcal{C}_{X}^{0,\alpha}(D)} &  = & \sup_{D}|u|+\sup
_{\overset{z,\zeta\in D}{z\neq\zeta}}\frac{|u(z)-u(\zeta)|}{d_{p}(z,\zeta)^{\alpha}%
},\\
\Vert u\Vert_{\mathcal{C}_{X}^{1,\alpha}(D)} &  = & \Vert u\Vert_{\mathcal{C}_{X}^{0,\%
alpha}(D)}+\sum_{i=1}^{q}\Vert X_{i}u\Vert_{\mathcal{C}_{X}^{0,\alpha}(D)}+
\sup_{\overset{z,\zeta\in D}{z\neq\zeta}}\frac{|u(z)-u(\zeta)-\sum_{j=1}^{q}%
(z_{j}-\zeta_{j})X_{j}u(\zeta)|}{d_{p}(z,\zeta)^{1+\alpha}},\\
\Vert u\Vert_{\mathcal{C}_{X}^{2,\alpha}(D)} &  = & \Vert u\Vert_{\mathcal{C}_{X}^{0,\alpha}(D)%
}+\sum_{i=1}^{q}\Vert X_{i}u\Vert_{\mathcal{C}_{X}^{0,\alpha}(D)}+\sum_{i,j=1}^{q}%
\Vert{X_{i}X_{j}}u\Vert_{\mathcal{C}_{X}^{0,\alpha}(D)%
}+\Vert\partial_{t}u\Vert_{\mathcal{C}_{X}^{0,\alpha}(D)}.
\end{eqnarray*}
In particular, note that we could just as well define H\"{o}lder spaces in terms of $\hat
{d}_{p}.$ Moreover, we let $\mathcal{C}^{0}(D)$ denote the set of functions
which are continuous on $D$. Note that if $u\in\mathcal{C}_{X}^{0,\alpha}(D)$,
then $u$ is H\"{o}lder continuous in the usual sense as we see using
\eqref{e-dist-eucl}. Let $m\in\{0,1,2\}$, $\alpha\in(0,1)$. If $u\in
\mathcal{C}_{X}^{k,\alpha}(D^{\prime})$ for every compact subset $D^{\prime}$
of $D$, then we write $u\in\mathcal{C}_{X,\text{loc}}^{k,\alpha}(D)$.
Furthermore, for $p\in\lbrack1,\infty]$
we define the Sobolev-Stein spaces
\[
\mathcal{S}_{X}^{p}(D)=\{u\in L^{p}(D):\ {X_{i}}u,\ {X_{i}X_{j}}%
u,\ \partial_{t}u\in L^{p}(D),\ i,j=1,...,q\}
\]
and we let
\[
\Vert u\Vert_{\mathcal{S}_{X}^{p}(D)}=\Vert u\Vert_{L^{p}(D)}+\sum_{i=1}%
^{q}\Vert{X_{i}}u\Vert_{L^{p}(D)}+\sum_{i,j=1}^{q}\Vert{X_{i}X_{j}}%
u\Vert_{L^{p}(D)}+\Vert\partial_{t}u\Vert_{L^{p}(D)}.
\]
Here $L^{p}(D)$ is the $L^{p}$-norm in $\mathbf{G}\equiv\mathbf{R}^{n+1}$ with
respect to the Lebesgue measure. If $u\in\mathcal{S}^{p}(D^{\prime})$ for
every compact subset $D^{\prime}$ of $D$, then we write $u\in\mathcal{S}%
_{X,loc}^{p}(D)$.

We next consider the spaces $\mathcal{C}_{X}^{m,\alpha}(D)$, for
$m\in\{0,1,2\}$, and we want to understand to which extent these spaces can be
characterized using polynomials. The following outline is based on the results
in \cite{GL} and \cite{B09}. Recall that $\delta_{\lambda}$ satisfies
\eqref{str1}. In the following we rewrite \eqref{str1} as
\begin{equation}
\delta_{\lambda}(x_{1},..,x_{n})=(\lambda^{\sigma_{1}}x_{1},..,\lambda
^{\sigma_{n}}x_{n})\label{str1+}%
\end{equation}
where $1\leq\sigma_{1}\leq...\leq\sigma_{n}<\infty$. Using this notation, yet
another example of a $\delta_{\lambda}$-homogeneous norm is given by
$||x||_{\mathbf{G}^{\prime}}=\sum_{j=1}^{n}|x_{j}|^{1/\sigma_{j}}$ whenever
$x\in\mathbf{R}^{n}$ . Furthermore, if $I=(i_{1},...,i_{n})$ is a multi-index,
$i_{j}\in\mathbf{N\cup}\{0\}$, we let
\[
|I|_{\delta_{\lambda}}=\sum_{j=1}^{n}\sigma_{j}i_{j}%
\]
denote the $\delta_{\lambda}$-length of $I$ and we denote, as usual, by
$|I|=\sum_{j=1}^{n}i_{j}$ the standard length of the multi-index $I$. A
polynomial $p(x)$ in the variable $x_{1},...,x_{n}$ is said to have
$\delta_{\lambda}$-degree $m$ if $p(x)=\sum_{|I|_{\delta_{\lambda}}\leq
m}a_{I}x^{I}$. We denote by $\mathbb{P}_{m,\delta_{\lambda}}$ the class of all
polynomials with $\delta_{\lambda}$-degree less than or equal to $m$. Let
$\Omega$ be any fixed bounded and open neighborhood of the origin and define
$\Omega_{\mu}=\delta_{\mu}\Omega$ for $\mu>0$. Given $1\leq\sigma_{1}%
\leq...\leq\sigma_{n}$ as in \eqref{str1+} we define $\Sigma=\{|I|_{\delta
_{\lambda}}:\ I=(i_{1},...,i_{n})\in(\mathbf{N}\cup\{0\})^{n},\ |I|_{\delta
_{\lambda}}>0\}$. Below $C^{m}(\mathbf{G},\mathbf{R})$ will denote the set of
functions $f:\mathbf{G\rightarrow R}$ for which differentials along vector
fields in $\mathbf{g}$ of $\delta_{\lambda}$-length $\leq m$ exists. We are
now ready to account for how to approximate functions using polynomials.
Below, $Z_{1},...,Z_{n}$ denotes the Jacobian basis of $\mathbf{g}$, the Lie
algebra of $\mathbf{G}$ and for $I=(i_{1},...,i_{n}),$ $Z^{I}$ denotes the
differential operator $Z^{I}=Z_{1}^{i_{1}}\cdots Z_{n}^{i_{n}}$. Similarly,
for $I=(i_{1},...,i_{k}),$ where $i_{j}\in\{1,...,n\}$ for $j=1,...,k,$ we
define the differential operator $Z_{I}=Z_{i_{1}}\cdots Z_{i_{k}}$. This time
we let $|I|_{\delta_{\lambda}}=\sum_{j=1}^{k}\sigma_{i_{j}}.$ Furthermore,
$B_{\mathbf{G}^{\prime}}(y,r)$ denotes $\{x\in\mathbf{R}^{n}:||y^{-1}\circ
x||_{\mathbf{G}^{\prime}}<r\}.$ For the proof of the following lemma we refer
to Theorem 2.8 in \cite{GL}.

\begin{lemma}
\label{polyy+} Let $\Omega\subset\mathbf{R}^{n}$ be open, let $u:\Omega
\rightarrow\mathbf{R}$ be a smooth function, $\alpha\in(0,1)$ and $m\in\Sigma
$. Suppose there exists positive constants $A$ and $\rho>0$ such that for each
$B_{\mathbf{G}^{\prime}}(z,R)\subset\Omega$ and for each $y\in B_{\mathbf{G}%
^{\prime}}(z,\rho R)$ there exists a polynomial $P_{m}^{y}\in\mathbb{P}%
_{m,\delta_{\lambda}}$ such that
\begin{equation}
|u(x)-P_{m}^{y}(x)|\leq A||y^{-1}\circ x||_{\mathbf{G}^{\prime}}^{m+\alpha
}\mbox{ for
all $x\in B_\mathbf{G \prime}(y,\rho R)$}.\label{star}%
\end{equation}
Then there exist positive constants $C$ and $\eta>0$ depending only on
$A,\alpha,m$,$\rho$ and the structure such that
\[
|Z^{I}u(x_{1})-Z^{I}u(x_{2})|\leq C||x_{1}^{-1}\circ x_{2}||_{\mathbf{G}%
^{\prime}}^{\alpha}%
\]
for all $|I|_{\delta_{\lambda}}=m$, $I=(i_{1},...,i_{n})$, and all
$x_{1},x_{2}\in B_{\mathbf{G}^{\prime}}(y,\eta R)$.
\end{lemma}

To find Taylor polynomials we use the following theorem, see Theorem 2 in
\cite{B09}.

\begin{theorem}
\label{Taylor}Let $\mathbf{G}=(\mathbf{R}^{n},\circ,\delta_{\lambda})$ be a
homogeneous Lie group. Let $y\in\mathbf{G}$ be fixed and assume that
$m\in\mathbf{N}\cup\{0\}$ and that $u\in C^{m+1}(\mathbf{G},\mathbf{R}).$ Let
$\{Z_{1},...,Z_{n}\}$ be the Jacobian basis for $\mathbf{g}$, the Lie algebra
of $\mathbf{G}.$ Then, for every $x\in\mathbf{G}$ the Taylor polynomial of
$\delta_{\lambda}$-degree $m$ of $u$ at $y$ is given by%
\begin{equation}
u(x)=u(y)+\sum\limits_{k=1}^{m}\underset{|I|_{\delta_{\lambda}}\leq m}%
{\sum_{I=(i_{1},...,i_{k})}}\frac{Z_{I}u(y)}{k!}\zeta_{i_{1}}(y_{{}}^{-1}\circ
x)\cdots\zeta_{i_{k}}(y^{-1}\circ x,t)+R_{m}(x,y),\label{Taylor1}%
\end{equation}
where the function $\zeta(h)$ is defined by%
\[
h=Exp(\zeta_{1}(h)Z_{1}+...+\zeta_{n}(h)Z_{n}).
\]
Moreover, there exists a constant $c$ depending only on $\mathbf{G}$ and
$d_{p}$ such that%
\begin{equation}
|R_{m}(x)|\leq\sum_{k=1}^{m+1}\frac{c^{k}}{k!}\underset{|I|_{\delta_{\lambda}%
}>m}{\sum_{I=(i_{1},...,i_{k})}}d_{p}(x,y)^{|I|_{\delta_{\lambda}}}\sup
_{||\xi||_{\mathbf{L}}\leq c||y^{-1}\circ x||_{\mathbf{L}}}|Z_{I}u(y\circ
\xi)|.\label{Taylor2}%
\end{equation}

\end{theorem}

\bigskip

Note that, in \cite{B09}, Bonfiglioli uses the distance $\hat{d}_{p}$ in
(\ref{Taylor2}). This theorem naturally extend to $\mathbf{L}$, and we remark
that $\zeta_{t}(h)=t$ while $\zeta_{i}(h)$ is independent of $t$ for
$i\in\{1,...,n\}$. We also mention that altough Lemma \ref{polyy+} and Theorem
\ref{Taylor} are stated using $Z^{I}$ and $Z_{I}$ respectively Proposition
20.1.4 in \cite{BLU07} shows, with slight abuse of notation, that $Z_{I}\in
span\{Z^{I}\}$. Given $(\xi,\tau)\in\mathbf{R}^{n+1}$ and for $(x,t)\in
\mathbf{R}^{n+1}$ we let
\begin{equation}
P_{0}^{(\xi,\tau)}u(x,t)\ \text{\ \ and \ }\ P_{1}^{(\xi,\tau)}%
u(x,t)\nonumber\label{po1}%
\end{equation}
denote the polynomials of $\delta_{\lambda}$-degree $m\in\{0,1\}$ such that
(\ref{star}) holds. Above $(\xi,\tau)$ is the point around which we
approximate $u$, i.e. it corresponds to $y$ in (\ref{star}). Using this
notation we prove:

\begin{lemma}
\label{polyy++} Let $m\in\{0,1\}$ and $\alpha\in(0,1)$. Then $u\in
\mathcal{C}_{X}^{m,\alpha}(D)$ if and only if there exists a constant $c$ such
that
\begin{equation}
\left\vert u(x,t)-P_{m}^{(\xi,\tau)}u(x,t)\right\vert \leq cd_{p}\left(
(x,t),(\xi,\tau)\right)  ^{m+\alpha}\label{Taylor3}%
\end{equation}
whenever $(x,t),\ (\xi,\tau)\in D$.
\end{lemma}

\noindent\textbf{Proof.} If $u\in\mathcal{C}_{X}^{m,\alpha}(D)$ then by
(\ref{Taylor1}) we have that%
\[
P_{0}^{(\xi,\tau)}u(x,t)=u(\xi,\tau)
\]
so (\ref{Taylor3}) holds with $c=||u||_{\mathcal{C}_{X}^{0,\alpha}(D)}$ by the
definition of $\mathcal{C}_{X}^{m,\alpha}(D).$ For $m=1$%
\[
P_{1}^{(\xi,\tau)}u(x,t)=u(\xi,\tau)+\sum_{i=1}^{q}X_{i}u(\xi,\tau)\xi_{i}%
((\xi,t)^{-1}\circ(x,t)),
\]
and since $\{X_{i}\}_{i=1}^{q}$ belong to the first layer of the
stratification and also to the Jacobian basis $\zeta_{i}((\xi,t)^{-1}%
\circ(x,t))=((\xi,t)^{-1}\circ(x,t))_{i}=x_{i}-\xi_{i}$ (note that this is not
the case for $i>q$). Once again (\ref{Taylor3}) holds with
$c=||u||_{\mathcal{C}_{X}^{1,\alpha}(D)}$ by the definition of $\mathcal{C}%
_{X}^{m,\alpha}(D).$

Now assume that (\ref{Taylor3}) holds. For $m=0$ (\ref{Taylor3}) reads%
\begin{equation}
|u(x,t)-u(\xi,\tau)|\leq cd_{p}((x,t),(\xi,\tau))^{\alpha},\label{ineq0}%
\end{equation}
and in particular $u\in\mathcal{C}_{X}^{0,\alpha}(D).$ For $m=1$
(\ref{Taylor3}) reads%
\begin{equation}
\left\vert u(x,t)-u(\xi,\tau)+\sum_{i=1}^{q}X_{i}u(\xi,\tau)(x_{i}-\xi
_{i})\right\vert \leq cd_{p}((x,t),(\xi,\tau))^{1+\alpha},\label{ineq1}%
\end{equation}
and by Lemma \ref{polyy+} $|X_{i}u(x,t)-X_{i}u(\xi,\tau)|\leq cd_{p}%
((x,t),(\xi,\tau))^{\alpha}.$ Using this in the inequality (\ref{ineq1}) we
find that there exists a constant $c,$ depending on $D$, such that
(\ref{ineq0}) holds and therefore it yields that $u\in\mathcal{C}%
_{X}^{1,\alpha}(D)$. \hfill$\Box$

\subsection{Function classes for the Cauchy-Dirichlet problem and the obstacle
problem\label{SecFC}}

In our proof of Theorem \ref{th-int} and Theorem \ref{th} we will make use of
certain estimates, at the initial state, for the Cauchy-Dirichlet problem
\begin{equation}%
\begin{cases}
\mathcal{H}u(x,t)=f(x,t), & \text{in}\ D,\\
u(x,t)=g(x,t), & \text{on}\ \partial_{p}D,
\end{cases}
\label{e-obsobs}%
\end{equation}
where $D\subset\mathbf{R}^{n+1}$ is a bounded cylindrical domain in the sense
of \eqref{cyl}. In the following we introduce certain function classes using
which we are able to give streamlined proofs of important lemmas and to prove
Theorem \ref{th-int} and Theorem \ref{th} in a unified way. Function classes
are introduced both for the Cauchy-Dirichlet problem and for the obstacle problem.

\begin{definition}
\label{d3.2x} Let $\mathcal{H}$ be defined as in \eqref{operator}, assume
\eqref{1.1x}-\eqref{1.2+-}, let $D\subset\mathbf{R}^{n+1}$ be a bounded
cylindrical domain in the sense of \eqref{cyl}, let $m\in\{0,1,2\}$,
$\alpha\in(0,1)$, and let $M_{1},M_{2},M_{3}$ be three positive constants.
Then we say that $(u,f,g)$ belongs to the class $\mathcal{D}_{m}%
(\mathcal{H},\alpha,D,M_{1},M_{2},M_{3})$ if $u$ is a solution to problem
\eqref{e-obsobs} with $f\in\mathcal{C}_{X}^{0,\alpha}(D)$, $g\in
\mathcal{C}_{X}^{m,\alpha}(\overline{D})$ and
\[
\Vert u\Vert_{L^{\infty}(D)}\leq M_{1},\quad\Vert f\Vert_{\mathcal{C}%
_{X}^{0,\alpha}(D)}\leq M_{2},\quad\Vert g\Vert_{\mathcal{C}_{X}^{m,\alpha
}(D)}\leq M_{3}.
\]

\end{definition}

\begin{definition}
\label{d3.2} Let $\mathcal{H}$ be defined as in \eqref{operator}, assume
\eqref{1.1x}-\eqref{1.2+-}, let $D\subset\mathbf{R}^{n+1}$ be a bounded
cylindrical domain in the sense of \eqref{cyl}, let $m\in\{0,1,2\}$,
$\alpha\in(0,1)$, and let $M_{1},M_{2},M_{3}$ be three positive constants. Let
$\varphi\in\mathcal{C}_{X}^{m,\alpha}(D)$, $g\in\mathcal{C}^{0}(\overline{D}%
)$, $g\geq\varphi$ on $\partial_{P}D$, and let $u$ be a strong solution to
problem \eqref{e-obs}. Then, for $m\in\{0,1,2\}$ we say that $(u,f,g,\varphi)$
belongs to the class $\mathcal{P}_{m}(\mathcal{H},\alpha,D,M_{1},M_{2},M_{3})$
if
\[
\Vert u\Vert_{L^{\infty}(D)}\leq M_{1},\quad\Vert f\Vert_{\mathcal{C}%
_{X}^{0,\alpha}(D)}\leq M_{2},\quad\Vert\varphi\Vert_{\mathcal{C}%
_{X}^{m,\alpha}(D)}\leq M_{3}.
\]

\end{definition}

\begin{definition}
\label{d3.2y} Let $\mathcal{H}$ be defined as in \eqref{operator}, assume
\eqref{1.1x}-\eqref{1.2+-}, let $D\subset\mathbf{R}^{n+1}$ be a bounded
cylindrical domain in the sense of \eqref{cyl}, let $m\in\{0,1,2\}$,
$\alpha\in(0,1)$, and let $M_{1},M_{2},M_{3},M_{4}$ be four positive
constants. Then, for $m\in\{0,1,2\}$ we say that $(u,f,g,\varphi)$ belongs to
the class $\mathcal{\tilde{P}}_{m}(\mathcal{H},\alpha,D,M_{1},M_{2}%
,M_{3},M_{4})$ if $u$ is a strong solution to problem \eqref{e-obs} with
$f\in\mathcal{C}_{X}^{0,\alpha}(D)$, $\varphi,g\in\mathcal{C}_{X}^{m,\alpha
}(\overline{D})$, $g\geq\varphi$ on $\partial_{P}D$ and
\[
\Vert u\Vert_{L^{\infty}(D)}\leq M_{1},\quad\Vert f\Vert_{\mathcal{C}%
_{X}^{0,\alpha}(D)}\leq M_{2},\quad\Vert g\Vert_{\mathcal{C}_{X}^{m,\alpha
}(D)}\leq M_{3},\quad\Vert\varphi\Vert_{\mathcal{C}_{X}^{m,\alpha}(D)}\leq
M_{4}.
\]

\end{definition}

As we will see below there are advantages, from the perspective of notation,
to introduce both of the classes $\mathcal{P}_{m}(\mathcal{H},\alpha
,D,M_{1},M_{2},M_{3})$ and $\mathcal{\tilde{P}}_{m}(\mathcal{H},\alpha
,D,M_{1},M_{2},M_{3},M_{4})$ though the classes are similar. In fact, we will
use the class $\mathcal{P}_{m}$ in the proof of Theorem \ref{th-int} and the
class $\mathcal{\tilde{P}}_{m}$ in the proof of Theorem \ref{th}. Furthermore,
we note that our proofs of Theorem \ref{th-int} and Theorem \ref{th} will be
based on certain blow-up arguments and in the following we lay out the
fundamental notation concerning translations and dilations. The reason this
can be done in an efficient way is due to the assumption concerning the
existence of the homogeneous Lie group $\mathbf{G}=(\mathbf{R}^{n}%
\!,\circ,\delta_{\lambda})$. In particular, let $D\subset\mathbf{R}^{n+1}$ be
a bounded cylindrical domain in the sense of \eqref{cyl}, assume that
$(0,0)\in D$, and let $v\in\mathcal{C}^{0}(D)$. We define the blow-up of a
function $v$, at $(0,0)$, as
\[
v^{r}(x,t)=v^{r,(0,0)}(x,t):=v\left(  \delta_{r}(x,t)\right)  ,\qquad r>0,
\]
whenever $\delta_{r}(x,t)\in D$. Recall that $\delta_{r}(x,t)$ was defined in
Section 2.1. Using this notation a direct computation shows that
\begin{equation}
\mathcal{H}u=f\ \text{in}\ D\quad\text{if and only if}\quad\mathcal{H}%
_{r}u^{r}=r^{2}f^{r}\ \text{in}\ \delta_{1/r}(D),\label{riscal}%
\end{equation}
where
\begin{equation}
\mathcal{H}_{r}=\sum_{i,j=1}^{q}a_{ij}^{r}{X_{i}^{r}X_{j}^{r}}+\sum_{i=1}%
^{q}rb_{i}^{r}{X_{i}^{r}}-\partial_{t},\label{3.4}%
\end{equation}
since the vector fields $\{X_{i}\}_{i=1}^{q}$ are $\delta_{\lambda}%
$-homogeneous of degree 1. In particular, $u^{r}(x,t)=u(\delta_{r}(x,t))$,
$a_{ij}^{r}(x,t)=a_{ij}(\delta_{r}(x,t))$, $b_{i}^{r}(x,t)=b_{i}(\delta
_{r}(x,t))$ whenever $\delta_{r}(x,t)\in D$. Furthermore, $X_{i}^{r}%
=(c_{i1}^{r}(x),...,c_{in}^{r}(x))\cdot(\partial_{x_{1}},...,\partial_{x_{n}%
})$, where $c_{ik}^{r}(x)=c_{ik}(\delta_{r}(x)) $. Note that $X^{r}%
=\{X_{1}^{r},...,X_{m}^{r}\}$ is still a system of smooth vector fields
satisfying \eqref{1.1x} and $A^{r}=\{a_{ij}^{r}\}$ and $b_{i}^{r}$ satisfy,
for $r\in(0,1]$, \eqref{1.2} and (\ref{1.2+-}) with the same constant as
$A=\{a_{ij}\}$ and $b_{i}$. To proceed we also define, given $r\in(0,1]$ and
$(x_{0},t_{0})\in\mathbf{R}^{n+1}$,
\begin{equation}
u^{r,(x_{0},t_{0})}(x,t)=u((x_{0},t_{0})\circ\delta_{r}(x,t)).\label{and-e1}%
\end{equation}
Let $m\in\{0,1,2\}$, $\alpha\in(0,1)$ and consider $r\in(0,1]$. We then remark
that $u\in\mathcal{C}_{X}^{m,\alpha}({D})$ if and only if $u^{r,(x_{0},t_{0}%
)}\in\mathcal{C}_{X}^{m,\alpha}(\delta_{1/r}(D))$ and we note that
\begin{equation}
\Vert u^{r,(x_{0},t_{0})}\Vert_{\mathcal{C}_{X}^{m,\alpha}(\delta_{1/r}%
((x_{0},t_{0})^{-1}\circ D))}\leq\Vert u\Vert_{\mathcal{C}_{X}^{m,\alpha}%
({D})}.\label{fredag}%
\end{equation}
Indeed in the case $m=0$ we have
\[
\Vert u^{r,(x_{0},t_{0})}\Vert_{\mathcal{C}_{X}^{m,\alpha}(\delta_{1/r}%
(D))}=\sup_{D}|u|+r^{\alpha}\sup_{\overset{z,\zeta\in D}{z\neq\zeta}}\frac
{|u(z)-u(\zeta)|}{d_{p}(z,\zeta)^{\alpha}}\leq\Vert u\Vert_{\mathcal{C}%
_{X}^{m,\alpha}({D})}.
\]
Moreover
\[
\mathcal{H}u=f\ \text{in}\ (x_{0},t_{0})\circ\delta_{r}(D)\quad\text{if and
only if}\quad\mathcal{H}_{r}^{(x_{0},t_{0})}u^{r,(x_{0},t_{0})}=r^{2}%
f^{r,(x_{0},t_{0})}\ \text{in}\ D,
\]
where
\[
\mathcal{H}_{r}^{(x_{0},t_{0})}=\sum_{i,j=1}^{q}a_{ij}^{r,(x_{0},t_{0}%
)}(x,t)X_{i}^{r,(x_{0},t_{0})}X_{j}^{r,(x_{0},t_{0})}+\sum_{i=1}^{q}%
rb_{i}^{r,(x_{0},t_{0})}(x,t)X_{i}^{r,(x_{0},t_{0})}-\frac{\partial}{\partial
t}.
\]
Above $a_{ij}^{r,(x_{0},t_{0})}$, $b_{i}^{r,(x_{0},t_{0})}$ and $X_{i}%
^{r,(x_{0},t_{0})}$ are defined as in (\ref{and-e1}). In particular, we can
conclude that if $(x_{0},t_{0})\in\mathbf{R}^{n+1}$, $r\in(0,1]$, and if
$(u,f,g,\varphi)\in\mathcal{P}_{m}(\mathcal{H},\alpha,D,M_{1},M_{2},M_{3})$,
then
\[
(u^{r,(x_{0},t_{0})},f^{r,(x_{0},t_{0})},g^{r,(x_{0},t_{0})},\varphi
^{r,(x_{0},t_{0})})\in\mathcal{P}_{m}(\mathcal{H}_{r}^{(x_{0},t_{0})}%
,\alpha,(x_{0},t_{0})\circ\delta_{r}(D),\alpha,M_{1},M_{2},M_{3}).
\]
The same statement holds for the class $\mathcal{\tilde{P}}_{m}(\mathcal{H}%
,\alpha,D,M_{1},M_{2},M_{3},M_{4})$.

\section{Estimates for Parabolic Non-divergence Operators of H\"{o}rmander
type\label{secFL}}

In this section we collect a number of results concerning\ parabolic
non-divergence operators of H\"{o}rmander type. These results will be used in
the proof of Theorem \ref{th-int} and Theorem \ref{th}. In particular, let
$\mathcal{H}$ be defined as in \eqref{operator} and assume
\eqref{1.1x}-\eqref{1.2+-}. We first state some results from \cite{BBLU09}. In
particular, in \cite{BBLU09} it is proved that there exists a fundamental
solution, $\Gamma$, for $\mathcal{H}$ on $\mathbf{R}^{n+1}$ with a number of
important properties: $\Gamma$ is a continuous function away from the diagonal
of $\mathbf{R}^{n+1}\times\mathbf{R}^{n+1}$ and $\Gamma(x,t,\xi,\tau)=0$ for
$t\leq\tau$. Moreover,
\[
\Gamma(\cdot,\cdot,\xi,\tau)\in\mathcal{C}_{X,loc}^{2,\alpha}(\mathbf{R}%
^{n+1}\setminus\{(\xi,\tau)\})\mbox{ for every fixed
$(\xi,\tau)\in\mathbf R^{n+1}$}
\]
and $\mathcal{H}(\Gamma(\cdot,\cdot,\xi,\tau))=0$ in $\mathbf{R}%
^{n+1}\setminus\{(\xi,\tau)\}$. For every $\psi\in C_{0}^{\infty}%
(\mathbf{R}^{n+1})$ the function
\[
w(x,t)=\int\limits_{\mathbf{R}^{n+1}}\Gamma(x,t,\xi,\tau)\psi(\xi,\tau)d\xi
d\tau
\]
belongs to $\mathcal{C}_{X,loc}^{2,\alpha}(\mathbf{R}^{n+1}\setminus
\{(\xi,\tau)\})$ and $\mathcal{H}w=\psi$ in $\mathbf{R}^{n+1}$. Furthermore,
let $\mu\geq0$ and $T_{2}>T_{1}$ be such that $(T_{2}-T_{1})\mu$ is small
enough, let $g\in\mathcal{C}_{X}^{0,\beta}(\mathbf{R}^{n}\times\lbrack
T_{1},T_{2}]),$ $0<\beta\leq\alpha$, and let $f\in\mathcal{C}^{0}%
(\mathbf{R}^{n})$ be such that $|g(x,t)|,|f(x)|\leq c\exp(\mu d(x,0)^{2})$ for
some constant $c>0$. Then the function
\[
u(x,t)=\int\limits_{\mathbf{R}^{n}}\Gamma(x,t,\xi,T_{1})f(\xi)d\xi
+\int\limits_{T_{1}}^{t}\int\limits_{\mathbf{R}^{n}}\Gamma(x,t,\xi,\tau
)g(\xi,\tau)d\xi d\tau,\ x\in\mathbf{R}^{n},\ t\in(T_{1},T_{2}],
\]
belongs to the class $\mathcal{C}_{X,loc}^{2,\alpha}(\mathbf{R}^{n}%
\times(T_{1},T_{2}))\cap\mathcal{C}^{0}(\mathbf{R}^{n}\times\lbrack
T_{1},T_{2}]) $. Moreover, $u$ solves the Cauchy-Dirichlet problem
\[
\mathcal{H}u=g\mbox{ in }\mathbf{R}^{n}\times(T_{1},T_{2}),\ u(\cdot
,T_{1})=f(\cdot)\mbox{ in }\mathbf{R}^{n}.
\]
Concerning the fundamental solution we will use the following uniform Gaussian
bounds, see Theorem 10.7 in \cite{BBLU09}.

\begin{lemma}
\label{Gaussbound} Let $\mathcal{H}$ be defined as in \eqref{operator} and
assume \eqref{1.1x}-\eqref{1.2+-}. Then the fundamental solution $\Gamma$ for
$\mathcal{H}$ on $\mathbf{R}^{n+1}$ satisfies the following estimates. There
exist a positive constant $C=C(\mathcal{H}),$ and for every $T>0,$ a positive
constant $c=c(\mathcal{H},T)$ such that if $0<t-\tau\leq T$ and $x,\xi
\in\mathbf{R}^{n}$, then
\begin{align*}
(i)  & c^{-1}|B(x,\sqrt{t-\tau})|^{-1}e^{-Cd(x,\xi)^{2}/(t-\tau)}\leq
\Gamma(x,t,\xi,\tau)\leq c|B(x,\sqrt{t-\tau})|^{-1}e^{-C^{-1}d(x,\xi
)^{2}/(t-\tau)},\\
(ii)  & |X_{i}\Gamma(\cdot,t,\xi,\tau)(x)|\leq c(t-\tau)^{-1/2}|B(x,\sqrt
{t-\tau})|^{-1}e^{-C^{-1}d(x,\xi)^{2}/(t-\tau)},\\
(iii)  & |X_{i}X_{j}\Gamma(\cdot,t,\xi,\tau)(x)|+|\partial_{t}\Gamma
(x,\cdot,\xi,\tau)(t)|\leq c(t-\tau)^{-1}|B(x,\sqrt{t-\tau})|^{-1}%
e^{-C^{-1}d(x,\xi)^{2}/(t-\tau)}.
\end{align*}

\end{lemma}

Let $Q$ be the homogeneous dimension of the underlying Lie group $\mathbf{G}$
as introduced in \eqref{e-Q}. Using $Q$, $(i)$ in Lemma \ref{Gaussbound} can
be rewritten as
\[
(i^{\prime})\ \ c^{-1}|{t-\tau}|^{-Q/2}e^{-Cd(x,\xi)^{2}/(t-\tau)}\leq
\Gamma(x,t,\xi,\tau)\leq c|{t-\tau}|^{-Q/2}e^{-C^{-1}d(x,\xi)^{2}/(t-\tau)},
\]
for $0<t-\tau\leq T.$ Recall the notion of cylinders in (\ref{pol3aa}); in a
similar fashion we now define
\begin{align*}
C_{r_{1},r_{2}}(x,t)  & =B_{d}(x,r_{1})\times(t-r_{2}^{2},t+r_{2}^{2}),\\
C_{r_{1},r_{2}}^{+}(x,t)  & =B_{d}(x,r_{1})\times(t,t+r_{2}^{2}),\\
C_{r_{1},r_{2}}^{-}(x,t)  & =B_{d}(x,r_{1})\times(t-r_{2}^{2},t),
\end{align*}
whenever $(x,t)\in\mathbf{R}^{n+1}$ and $r_{1},r_{2}>0$. The following Harnack
inequality is proved in Theorem 15.1 in \cite{BBLU09}.

\begin{lemma}
\label{Harnack} Let $\mathcal{H}$ be defined as in \eqref{operator} and assume
\eqref{1.1x}-\eqref{1.2+-}. Let $R>0$, $0<h_{1}<h_{2}<1$ and $\gamma\in(0,1)$.
Then there exists a positive constant $c=c(h_{1},h_{2},\gamma,R) $ such that
the following holds for every $(\xi,\tau)\in\mathbf{R}^{n+1}$, $r\in(0,R]$.
If
\[
u\in\mathcal{C}_{X}^{2,\alpha}(C_{r}^{-}(\xi,\tau))\cap\mathcal{C}%
^{0}(\overline{C_{r}^{-}(\xi,\tau)})
\]
satisfies $\mathcal{H}u=0$, $u\geq0$, in $C_{r}^{-}(\xi,\tau)$, then
\[
u(x,t)\leq cu(\xi,\tau)\ \text{whenever }(x,t)\in\overline{B_{d}(\xi,\gamma
r)}\times\left[  \tau-h_{2}r^{2},\tau-h_{1}r^{2}\right]  .
\]

\end{lemma}

In addition we will need the following Schauder estimate for the operator
$\mathcal{H}$, which is proved in Theorem 1.1 in \cite{BB07}.

\begin{lemma}
\label{t-schauder} Let $\mathcal{H}$ be defined as in \eqref{operator} and
assume \eqref{1.1x}-\eqref{1.2+-}. Let $R>0$ and $(x,t)\in\mathbf{R}^{n+1}$.
If $u\in\mathcal{C}_{X,loc}^{2,\alpha}(C_{R}(x,t))$ satisfies $\mathcal{H}u=f$
in $C_{R}(x,t)$, for some function $f\in\mathcal{C}_{X}^{0,\alpha}%
(C_{R}(x,t))$, then there exists a positive constant $c$, depending on
$\mathcal{H}$, $\alpha$ and $R$ such that
\[
\Vert u\Vert_{\mathcal{C}_{X}^{2,\alpha}(C_{R/2}(x,t))}\leq c\left(  \Vert
u\Vert_{L^{\infty}(C_{R}(x,t))}+\Vert f\Vert_{\mathcal{C}_{X}^{0,\alpha}%
(C_{R}(x,t))}\right)  .
\]

\end{lemma}

We will also need the weak maximum principle for the operator $\mathcal{H}$
and we will use one due to Picone, e.g., see Theorem 13.1 in \cite{BBLU09}.

\begin{lemma}
\label{weakmax} Let $\mathcal{H}$ be defined as in \eqref{operator} and assume
\eqref{1.1x}-\eqref{1.2+-}. Let $D\subset\mathbf{R}^{n+1}$ be a bounded
cylindrical domain in the sense of \eqref{cyl}. Then the following holds; if
$u\in\mathcal{C}_{X}^{2,\alpha}(D)$, $\mathcal{H}u\geq0$ in $D$ and $\limsup
u\leq0$ on $\partial_{p}D$, then $u\leq0$ in $D$.
\end{lemma}

Finally, we have to ensure that the continuous Cauchy-Dirichlet problem is
solvable in the cylinders $C_{r}(x,t)$, $C_{r}^{+}(x,t)$, $C_{r}^{-}(x,t)$
introduced in \eqref{pol3aa}. However, as previously mentioned, this is not
the case in general and we will have to work with modified cylinders. To
explain this further, in \cite{U}, Uguzzoni develops what he refers to as a
"cone criterion" for the solvability of the Cauchy-Dirichlet problem for
non-divergence operators structured on H\"{o}rmander vector fields. This is a
generalization of the well-known positive density condition in classical
potential theory. In the following we describe his results in the setting of
cylindrical domains $D=\Omega\times(T_{1},T_{2}),$ where $\Omega$ is an open,
bounded domain in $\mathbf{R}^{n}.$ A bounded, open set $\Omega$ is said to
have outer positive $d$-density at $x_{0}\in\partial\Omega$ if there exist
$r_{0},\ \theta>0$ such that
\[
|B_{d}(x_{0},r)\backslash\overline{\Omega}|\geq\theta|B_{d}(x_{0}%
,r)|,\ \ \ \text{for all }r\in(0,r_{0}),
\]
where $|\cdot|$ denotes Euclidean volume. Furthermore, if $\Omega$ has outer
positive $d$-density at all points $x\in\partial\Omega$ and if $r_{0}$ and
$\theta$ can be chosen independent of $x,$ then $\Omega$ is said to satisfy
the uniform outer positive $d$-density condition. The following lemma is a
special case of Lemma 4.1 in \cite{U}.

\begin{lemma}
\label{Dirichlet0} Let $\mathcal{H}$ be defined as in \eqref{operator}, assume
\eqref{1.1x}-\eqref{1.2+-} and that $f$ and $g$ are continuous and bounded.
Let $D\subset\mathbf{R}^{n+1}$ be a bounded cylindrical domain of the specific
form $D=\Omega\times(T_{1},T_{2}).$ Assume that $\Omega$ satisfies the uniform
outer positive $d$-density condition. Then there exists a unique solution
$u\in\mathcal{C}_{X}^{2,\alpha}(D)\cap\mathcal{C}^{0}(\overline{D})$ to the
problem
\[
\mathcal{H}{u}=f\mbox{ in $D$, $u=g$ on $\partial_p D$},
\]
where $\alpha$ is the H\"{o}lder exponent in (\ref{1.2+-}).
\end{lemma}

\noindent The above lemma only applies for bounded domains satisfying the
outer positive $d$-density condition and in general the cylinders $C_{r}(x,t)
$, $C_{r}^{+}(x,t)$, $C_{r}^{-}(x,t)$ do not satisfy this criterion. However,
in Theorem 6.5 in \cite{LU}, Lanconelli and Uguzzoni prove that we can
approximate these cylinders with modified versions which are regular for the
Cauchy-Dirichlet problem.

\begin{theorem}
\label{Solvable} Let $\Omega$ be a bounded domain in $\mathbf{R}^{n}.$ Then,
for every $\delta>0$ there exist $\Omega_{\delta}$ such that $\{x\in
\Omega:d(x,\partial\Omega)>\delta\}\subset\Omega_{\delta}\subset\Omega$ and
$\Omega_{\delta}$ satisfies the uniform outer positive $d$-density condition.
\end{theorem}

\begin{remark}
\label{remark1}\noindent Below, we will assume that we work with modified
cylinders which are regular for the Cauchy-Dirichlet problem.
\end{remark}

\hfill

\subsection{Auxiliary technical estimates}

In this subsection we derive a few auxiliary technical estimates to be used in
the forthcoming sections.

\begin{lemma}
\label{cor} Let $\mathcal{H}$ be defined as in \eqref{operator} with
corresponding fundamental solution $\Gamma$ and assume
\eqref{1.1x}-\eqref{1.2+-}. Given $\gamma,R>0$, we define the function
\[
u(x,t)=\int\limits_{B_{d}(0,R)}\Gamma(x,t,y,0)||y||^{\gamma}dy,\qquad\text{\ for
}x\in\mathbf{R}^{n},\ t>0.
\]
Let $D=\Omega\times(0,T)$ where $\Omega\subset\mathbf{R}^{n}$ is a bounded
domain and let $T>0$. Then there exists a positive constant $c=c(\gamma
,T,\mathcal{H})$ such that
\[
u(x,t)\leq c\Vert(x,t)\Vert_{p}^{\gamma},\qquad\text{for }(x,t)\in D.
\]

\end{lemma}

\noindent\textbf{Proof.} First of all we note that%
\[
||y||=||(y,0)||_{p}=||(y,0)^{-1}||_{p}\leq||(y,0)^{-1}\circ(x,t)||_{p}%
+||(x,t)||_{p}%
\]
by the triangle inequality. Hence,%
\begin{align}
u(x,t)  & =\int\limits_{B_{d}(0,R)}\Gamma(x,t,y,0)||y||^{\gamma}dy\nonumber\\
& \leq c(\gamma)\left(  ||(x,t)||_{p}^{\gamma}\int\limits_{B_{d}%
(0,R)}\Gamma(x,t,y,0)dy+\int\limits_{B_{d}(0,R)}\Gamma(x,t,y,0)||(y,0)^{-1}%
\circ(x,t)||_{p}^{\gamma}dy\right)  .\label{uxt1}%
\end{align}
Below we will use Lemma \ref{Gaussbound} and the following estimates, see
Proposition 10.11 and Corollary 10.12 in \cite{BBLU09}: For any $\mu\geq0$
there exist $c(\mu)$ such that%
\[
(d(x,y)^{2}/t)^{\mu}|B(x,\sqrt{\lambda t})|^{-1}\exp\left(  -\frac{d(x,y)^{2}%
}{\lambda t}\right)  \leq c(\mu,\mathcal{H})\lambda^{\mu}|B(x,\sqrt{2\lambda
t})|^{-1}\exp\left(  -\frac{d(x,y)^{2}}{2\lambda t}\right)
\]
for every $\lambda>0,$ $x,y\in\mathbf{R}^{n},\ t>0,$ and%
\[
\int_{\mathbf{R}^{n}}|B(x,\sqrt{ct})|^{-1}\exp\left(  -\frac{d(x,y)^{2}}%
{ct}\right)  dy\leq c(T,\mathcal{H})
\]
for $0<t\leq T,\ x\in\mathbf{R}^{n}$. Continuing, using the above estimates,
the first integral in (\ref{uxt1}) is bounded by some constant depending on
$T$ while the second integral is bounded by%
\begin{align*}
& \int\limits_{B_{d}(0,R)}\Gamma(x,t,y,0)||(y,0)^{-1}\circ(x,t)||_{p}^{\gamma
}dy\\
& \leq\int\limits_{B_{d}(0,R)}|B(x,\sqrt{\mathbf{c}t})|^{-1}\exp\left(
-\frac{d(x,y)^{2}}{\mathbf{c}t}\right)  \left(  d(x,y)^{2}+t\right)
^{\gamma/2}dy\\
& \leq c(\gamma,\mathcal{H})t^{\gamma/2}\left(  \int\limits_{B_{d}%
(0,R)}|B(x,\sqrt{\mathbf{c}t})|^{-1}\exp\left(  -\frac{d(x,y)^{2}}%
{\mathbf{c}t}\right)  \left(  d(x,y)^{2}/t\right)  ^{\gamma/2}dy+\right. \\
&
\ \ \ \ \ \ \ \ \ \ \ \ \ \ \ \ \ \ \ \ \ \ \ \ \ \ \ \ \ \ \ \ \ \ \ \ \ \ \ \ \ \ \ \ \ \ \left.
\int\limits_{B_{d}(0,R)}|B(x,\sqrt{\mathbf{c}t})|^{-1}\exp\left(
-\frac{d(x,y)^{2}}{\mathbf{c}t}\right)  dy\right) \\
& \leq c(\gamma,T,\mathcal{H})t^{\gamma/2}\leq c(\gamma,T,\mathcal{H}%
)||(x,t)||_{p}^{\gamma}.
\end{align*}
\hfill\noindent\noindent The dependency on $\mathcal{H}$ stems from the
constants in Lemma \ref{Gaussbound}. This concludes the proof.\hfill$\Box
$\newline

\begin{lemma}
\label{t-2.19} Let $\mathcal{H}$ be defined as in \eqref{operator} and assume
\eqref{1.1x}-\eqref{1.2+-}. Let $Q$ be the homogeneous dimension of the
underlying Lie group $\mathbf{G}$ as introduced in \eqref{e-Q}. Let $K\gg1$
and $R$ be given. Assume that $\mathcal{H}u=0$ in $C_{KR,R}^{-}(0,0)$ and that
$u(x,-R^{2})=0$ whenever $x\in\overline{B_{d}(0,KR)}$. Then there exists a
constant $c=c(\mathcal{H},T)$ such that
\[
\sup_{C_{R}^{-}}|u|\leq cK^{Q+1}e^{-cK^{2}}\sup_{\partial_{p}C_{KR,R}%
^{-}(0,0)\cap\{(x,t):\ t>-R^{2}\}}|u|,
\]

\end{lemma}

\noindent\textbf{Proof.} Let $r>0$ be small enough so that
$\{y:||y||<2r\}\subset B_{d}(0,1)$ and let $\phi\in C^{\infty}(\mathbf{R}%
^{n})$ be a function taking values in $[0,1]$ such that%
\[
\left\{
\begin{array}
[c]{cc}%
\phi(x)=1, & \text{if }||x||\geq2r,\\
\phi(x)=0, & \text{if }||x||\leq r.
\end{array}
\right.
\]
We define
\[
\omega(x,t):=2\int_{\mathbf{R}^{n}}\Gamma(x,t,y,-R^{2})\phi(\delta
_{1/KR}(y))dy,
\]
and we note that $\omega$ is a non-negative solution to the Cauchy-problem%
\[
\left\{
\begin{array}
[c]{ll}%
\mathcal{H}u=0, & \text{in }\mathbf{R}^{n}\times(-R^{2},0],\\
u(x,-R^{2})=2\phi(\delta_{1/KR}(x)), & \text{on }\mathbf{R}^{n}.
\end{array}
\right.
\]
Furthermore, we note that for $(x,t)\in\partial_{p}C_{KR}^{-}(0,0)\cap
\{(x,t):t>-R^{2}\}$ the dilated point $\delta_{1/KR}(x,t)$ belongs to
$C_{1,1/K^{2}}^{-}(0,0)$ and that, by Lemma \ref{Gaussbound},
\[
\omega(x,t)\geq2c^{-1}\int_{\mathbf{R}^{n}}\left\vert B\left(  x,\sqrt
{t+R^{2}}\right)  \right\vert ^{-1}\exp\left(  -C\frac{d(x,y)^{2}}{t+R^{2}%
}\right)  \phi(\delta_{1/KR}(y))dy.
\]
We intend to prove that the integral above converges to $2$ as $K\rightarrow
\infty$ uniformly in $(x,t)\in\partial_{p}C_{KR}^{-}(0,0)\cap\{(x,t):t>-R^{2}%
\}.$ Since, in that case, there exists a constant $K_{0}$ such that, for
$K>K_{0},$ we have that $\omega(x,t)\geq1$ for all $(x,t)\in\partial_{p}%
C_{KR}^{-}(0,0)\cap\{(x,t):t>-R^{2}\}.$ Then using the assumptions the weak
maximum principle implies that%
\[
u(x,t)\leq\omega(x,t)\sup_{\partial_{p}C_{KR,R}^{-}(0,0)\cap\{(x,t):\ t>-R^{2}%
\}}|u|,
\]
for boundary points and we are done if we, in addition, can prove that
$\omega(x,t)\leq cK^{Q+1}e^{-cK^{2}}$ for $(x,t)\in C_{R}^{-}(0,0)\cap
\{(x,t):t>-R^{2}\}.$

Now,
\begin{align*}
|\omega(x,t)-1|  & \leq c^{-1}\int_{\mathbf{R}^{n}}\left\vert B\left(
x,\sqrt{t+R^{2}}\right)  \right\vert ^{-1}\exp\left(  -C\frac{d(x,y)^{2}%
}{t+R^{2}}\right)  |\phi(\delta_{1/KR}(y))-1|dy\\
& \leq c^{-1}\int_{B_{d}(0,2rKR)}\left\vert B\left(  x,\sqrt{t+R^{2}}\right)
\right\vert ^{-1}\exp\left(  -C\frac{d(x,y)^{2}}{t+R^{2}}\right)  dy\\
& \leq c^{-1}R\int_{B_{d}(0,2rKR)}\left(  \sqrt{C}\frac{d(x,y)}{\sqrt{t+R^{2}%
}}\right)  ^{Q+1}\exp\left(  -C\frac{d(x,y)^{2}}{t+R^{2}}\right)  \left(
\frac{1}{\sqrt{C}d(x,y)}\right)  ^{Q+1}dy\\
& \leq c^{-1}Rc(Q,\mathcal{H})\int_{B_{d}(0,2rKR)}d(x,y)^{-Q-1}dy\\
& \leq c^{-1}Rc(Q,\mathcal{H})(2rKR)^{Q}\frac{1}{\left(  (\frac{1}{c_{d}%
}-2r)KR\right)  ^{Q+1}}\\
& =c(\Omega,Q_{T},\mathcal{H)}\frac{1}{K}\rightarrow0
\end{align*}
as $K\rightarrow\infty$ for $(x,t)\in C_{KR,R}^{-}(0,0)\cap\{(x,t):t>-R^{2}%
\}.$ Above we used that $2r<1.$ It remains to show that $\omega(x,t)\leq
cK^{Q+1}e^{-cK^{2}}$ for $(x,t)\in C_{R}^{-}(0,0)\cap\{(x,t):t>-R^{2}\}.$ To
prove this, recall that $u^{n}e^{-u^{2}}$ is bounded. A direct calculation
then shows that
\begin{align}
\omega(x,t)  & \leq c\int_{\mathbf{||y||}\geq rKR}\left\vert B\left(
x,\sqrt{t+R^{2}}\right)  \right\vert ^{-1}\exp\left(  -C\frac{d(x,y)^{2}%
}{t+R^{2}}\right)  dy\nonumber\\
& \leq c\int_{\mathbf{||y||}\geq rKR}(t+R^{2})^{-(Q+1)/2}\exp\left(  -\frac
{C}{4}\frac{d(0,y)^{2}}{(t+R^{2})}\right)  dy\nonumber\\
& \leq c\int_{\mathbf{||y||}\geq rKR}\left(  \frac{C}{8}\frac{d(0,y)^{2}%
}{K^{2}(t+R^{2})}\right)  ^{(Q+1)/2}\exp\left(  -\frac{C}{8}\frac{d(0,y)^{2}%
}{K^{2}(t+R^{2})}\right)  \cdot\nonumber\\
&
\ \ \ \ \ \ \ \ \ \ \ \ \ \ \ \ \ \ \ \ \ \ \ \ \ \ \ \ \ \ \ \ \ \ \ \ \ \ \ \ \ \ \ \ \cdot
\ \exp\left(  -\frac{C}{8}\frac{d(0,y)^{2}}{t+R^{2}}\right)  \left(  \frac
{8}{C}\frac{K^{2}}{d(0,y)^{2}}\right)  ^{(Q+1)/2}dy\nonumber\\
& \leq cK^{Q+1}\int_{\mathbf{||y||}\geq rKR}d(0,y)^{-(Q+1)}\exp\left(
-\frac{C}{8}\frac{d(0,y)^{2}}{t+R^{2}}\right)  dy\nonumber\\
& \leq cK^{Q+1}e^{-CK^{2}}\int_{\mathbf{||y||}\geq rKR}d(0,y)^{-(Q+1)}%
dy.\label{pet}%
\end{align}
Now, we define $S_{j}:=\{y:2^{j}rKR\leq||y||\leq2^{j+1}rKR\}$ for
$j\in\{0,1,...\}$, and we use (\ref{pet}) together with (\ref{bvolym}), to
obtain
\begin{align*}
\omega(x,t)  & \leq c(\mathcal{H},T\mathcal{)}K^{Q+1}e^{-CK^{2}}\sum
_{j=0}^{\infty}(2^{j}rKR)^{Q}(2^{j}rKR)^{-(Q+1)}\\
& \leq c(\mathcal{H},T\mathcal{)}K^{Q+1}e^{-CK^{2}}\sum_{j=0}^{\infty}%
2^{-j}\leq C(\mathcal{H},R_{0}\mathcal{)}K^{Q+1}e^{-CK^{2}},
\end{align*}
which completes the proof of the lemma. It should be clear that the constant
only depend on $\mathcal{H}$ and $T$, and that this dependency originates in
the use of Lemma \ref{Gaussbound}.\hfill$\Box$\newline

\begin{lemma}
\label{lemD} Let $\mathcal{H}$ be defined as in \eqref{operator} and assume
\eqref{1.1x}-\eqref{1.2+-}. Let $\alpha\in(0,1)$, $m=0,1,2$ and let
$M_{1},M_{2},M_{3}$ be positive constants. Assume that
\[
(u,f,g)\in{\mathcal{D}}_{m}(\mathcal{H},\alpha,C^{+},M_{1},M_{2},M_{3}).
\]
Then there exists a constant $c=c(\mathcal{H},m,\alpha,M_{1},M_{2},M_{3})$
such that
\[
\sup_{C_{r}^{+}}|u-g|\leq cr^{\gamma},\ r\in(0,1),
\]
where $\gamma=m+\alpha$ for $m=0,1$ and $\gamma=2$ for $m=2$.
\end{lemma}

\noindent\textbf{Proof.} Recall that $C^{+}=C_{1}^{+}(0,0)$ and that
$C_{r}^{+}=C_{r}^{+}(0,0)$. To start the proof we first note, using the
triangle inequality and Lemma \ref{polyy++}, that it suffices to prove that
\[
\sup_{C_{r}^{+}}\left\vert u-P_{m}^{(0,0)}g\right\vert \leq cr^{\gamma},\qquad
r\in(0,1),
\]
where $\gamma=m+\alpha$ for $m=0,1$ and that
\[
\sup_{C_{r}^{+}}\left\vert u-g\right\vert \leq cr^{\gamma},\qquad r\in(0,1),
\]
where $\gamma=2$ for $m=2$. Consider $m=0,1$ and let $v_{m}=u-P_{m}^{(0,0)}g
$. Then $v_{m}$ satisfies the equation
\[
\mathcal{H}v_{m}=f-\mathcal{H}P_{m}^{(0,0)}g=:f_{m}.
\]
Since $f_{m}\in\mathcal{C}_{X}^{0,\alpha}$ we see that we can assume in case
$m=0,1$, without loss of generality, that $P_{m}^{(0,0)}g=0$. By subtraction
of $g$ we can, using the same argument, assume that $g\equiv0$ in case $m=2$.
After these preliminary problem reductions we proceed with the proof and we
first consider the case $m=0$. Let
\[
\partial_{p}^{+}C^{+}=\partial_{p}C^{+}\cap\{t>0\},\qquad\partial_{p}^{-}%
C^{+}=\partial_{p}C^{+}\cap\{t=0\}.
\]
Let $v_{1},v_{2}$ and $v_{3}$ respectively be the unique solutions of the
following boundary value problems,
\[%
\begin{cases}
\mathcal{H}v_{1}=0\  & \text{in}\ C^{+},\\
v_{1}=0\  & \text{on}\ \partial_{p}^{+}C^{+},\\
v_{1}=g\  & \text{on}\ \partial_{p}^{-}C^{+},
\end{cases}
\qquad%
\begin{cases}
\mathcal{H}v_{2}=0\  & \text{in}\ C^{+},\\
v_{2}=g\  & \text{on}\ \partial_{p}^{+}C^{+},\\
v_{2}=0\  & \text{on}\ \partial_{p}^{-}C^{+},
\end{cases}
\qquad%
\begin{cases}
\mathcal{H}v_{3}=-\Vert f\Vert_{\mathcal{C}_{X}^{0,\alpha}(C^{+})}\  &
\text{in}\ C^{+},\\
v_{3}=0\  & \text{on}\ \partial_{p}C^{+}.
\end{cases}
\]
Note that using Lemma \ref{Dirichlet0}, Theorem \ref{Solvable} and Remark
\ref{remark1} we see that $C^{+}$ is regular for the Cauchy-Dirichlet problem
for $\mathcal{H}$ and that $v_{1},v_{2}$ and $v_{3}$ are well-defined. Then,
by the maximum principle, see Lemma \ref{weakmax}, we have
\[
v_{1}+v_{2}-v_{3}\leq u\leq v_{1}+v_{2}+v_{3}\qquad\text{in }C^{+},
\]
so that we only have to prove, since $\gamma=\alpha$ in the case $m=0$, that
\begin{equation}
\sup_{C_{r}^{+}}\left(  \left\vert v_{1}\right\vert +\left\vert v_{2}%
\right\vert +\left\vert v_{3}\right\vert \right)  \leq cr^{\alpha
},\label{3.41b}%
\end{equation}
for $r$ suitably small. To proceed, since $\Vert g\Vert_{\mathcal{C}%
_{X}^{0,\alpha}(C^{+})}\leq M_{3}$ we see by the maximum principle that
\[
|v_{1}(x,t)|\leq M_{3}\int\limits_{\mathbf{R}^{n}}\Gamma(x,t,y,0)\Vert
y\Vert^{\alpha}dy
\]
and hence we can apply Lemma \ref{cor} to conclude that
\[
|v_{1}(x,t)|\leq cM_{3}\Vert(x,t)\Vert_{p}^{\alpha}.
\]
We next apply Lemma \ref{t-2.19} to conclude that
\[
\sup_{C_{r}^{+}}|v_{2}|\leq cr^{Q}\exp\left(  -{c}r^{-2}\right)
\sup_{\partial_{p}^{+}C^{+}}|v_{2}|
\]
for any $r\leq K^{-1}$. Since $|v_{2}|$ agrees with $|u|$ on $\partial_{p}%
^{+}C^{+}$ we can conclude that
\[
\sup_{C_{r}^{+}}|v_{2}|\leq cM_{1}r^{Q}\exp\left(  -{c_{1}}r^{-2}\right)  \leq
c_{2}M_{1}r^{2},\qquad\text{for every}\quad r\in(0,K^{-1}).
\]
Finally, we have
\[
|v_{3}(x,t)|\leq\Vert f\Vert_{\mathcal{C}_{X}^{0,\alpha}(C^{+})}%
\int\limits_{0}^{t}\int\limits_{\mathbf{R}^{n}}\Gamma(x,t,y,s)dyds\leq
ct\,\Vert f\Vert_{\mathcal{C}_{X}^{0,\alpha}(C^{+})}\leq cM_{2}\Vert
(x,t)\Vert_{p}^{2}.
\]
This proves \eqref{3.41b} and the proof of the lemma is complete in the case
$m=0$. The cases $m=1,2$ can be handled by repeating the argument above. In
particular, we now apply Lemma \ref{cor} with $\gamma=m+\alpha$ for $m=1$ and
$\gamma=2$ for $m=2$ and we find
\[
\left\vert v_{1}(x,t)\right\vert \leq cM_{3}\Vert(x,t)\Vert_{p}^{\gamma}.
\]
Using this estimate we are then able to complete the argument as above. This
completes the proof.\hfill$\Box$

\section{Estimates for the obstacle problem}

The purpose of this section is to derive certain estimates for the obstacle
problem to be used in the proofs of Theorem \ref{th-int} and Theorem \ref{th}.
The results of the section are the following lemmas.

\begin{lemma}
\label{lem1} Let $\mathcal{H}$ be defined as in \eqref{operator} and assume
\eqref{1.1x}-\eqref{1.2+-}. Let $\alpha\in(0,1)$, $m=0,1,2$, and let
$M_{1},M_{2},M_{3}$ be positive constants. Assume that
\[
(u,f,g,\varphi)\in{\mathcal{P}}_{m}(\mathcal{H},\alpha,C^{-},M_{1},M_{2}%
,M_{3})\mbox{ and that $u(0,0)=\varphi(0,0)$}.
\]
Then there exists $c=c(\mathcal{H},m,\alpha,M_{1},M_{2},M_{3})$ such that
\[
\sup_{C_{r}^{-}}|u-\varphi|\leq cr^{\gamma},\qquad r\in(0,1),
\]
where $\gamma=m+\alpha$ for $m=0,1$ and $\gamma=2$ for $m=2$.
\end{lemma}

\begin{lemma}
\label{lem} Let $\mathcal{H}$ be defined as in \eqref{operator} and assume
\eqref{1.1x}-\eqref{1.2+-}. Let $\alpha\in(0,1)$, $m=0,1,2$, and let
$M_{1},M_{2},M_{3},M_{4}$ be positive constants. Assume that
\[
(u,f,g,\varphi)\in{\mathcal{\tilde{P}}}_{m}(\mathcal{H},\alpha,C^{+}%
,M_{1},M_{2},M_{3},M_{4}).
\]
Then there exists $c=c(\mathcal{H},m,\alpha,M_{1},M_{2},M_{3},M_{4})$ such
that
\[
\sup_{C_{r}^{+}}|u-g|\leq cr^{\gamma},\qquad r\in(0,1),
\]
where $\gamma=m+\alpha$ for $m=0,1$ and $\gamma=2$ for $m=2$.
\end{lemma}

\begin{lemma}
\label{lem1+} Let $\mathcal{H}$ be defined as in \eqref{operator} and assume
\eqref{1.1x}-\eqref{1.2+-}. Let $\alpha\in(0,1)$, $m=0,1,2$, and let
$M_{1},M_{2},M_{3}$ be positive constants. Assume that
\[
(u,f,g,\varphi)\in{\mathcal{P}}_{m}(\mathcal{H},\alpha,C,M_{1},M_{2}%
,M_{3})\mbox{
and that $u(0,0)=\varphi(0,0)$}.
\]
Then there exists $c=c(\mathcal{H},m,\alpha,M_{1},M_{2},M_{3})$ such that
\[
\sup_{C_{r}}|u-\varphi|\leq cr^{\gamma},\qquad r\in(0,1),
\]
where $\gamma=m+\alpha$ for $m=0,1$ and $\gamma=2$ for $m=2$.
\end{lemma}

\subsection{Proof of Lemma \ref{lem1}}

We first note, as in the proof of Lemma \ref{lemD}, that using the triangle
inequality and Lemma \ref{polyy++} it suffices to prove that
\begin{equation}
\sup_{C_{r}^{-}}\left\vert u-P_{m}^{(0,0)}\varphi\right\vert \leq
cr^{\gamma},\qquad r\in(0,1),\label{a}%
\end{equation}
where $\gamma=m+\alpha$ for $m=0,1$ and
\begin{equation}
\sup_{C_{r}^{-}}\left\vert u-\varphi\right\vert \leq cr^{\gamma},\qquad
r\in(0,1),\label{b}%
\end{equation}
where $\gamma=2$ for $m=2$. Furthermore, by a simple problem reduction we see
that we can assume, without loss of generality, that $P_{m}^{(0,0)}%
\varphi\equiv0$ for $m=0,1$ and that $\varphi\equiv0$ for $m=2$. After these
simplifications we define $S_{k}^{-}(u)$ as in \eqref{1.17x} and we intend to
prove that there exists a positive constant $\tilde{c}=\tilde{c}\left(
\mathcal{H},m,\alpha,M_{1},M_{2},M_{3}\right)  $ such that \eqref{1.17z} holds
for all $k\in\mathbf{N}$ with $F$ defined as in \eqref{lex1}. Indeed, if
\eqref{1.17z} holds then we see, by a simple iteration argument, that
\[
S_{k}^{-}(u-F)\leq\frac{\tilde{c}}{2^{k\gamma}}%
\]
and hence \eqref{a} and \eqref{b}, and thereby Lemma \ref{lem1}, follows.

We first consider the case when\emph{\ }$m=0$\emph{\ }and prove \eqref{1.17z}
with\emph{\ }$\gamma=\alpha$\emph{.} In particular, we assume that
\[
(u,f,g,\varphi)\in{\mathcal{P}}_{0}(\mathcal{H},\alpha,C^{-},M_{1},M_{2}%
,M_{3}).
\]
As in \cite{FNPP} we divide the argument into three steps.

\medskip\noindent\textbf{Step 1.} Using that $\varphi(0,0)=0$ by the problem
reduction, we first note that
\begin{equation}
u(x,t)\geq\varphi(x,t)=\varphi(x,t)-\varphi(0,0)\geq-M_{3}\Vert(x,t)\Vert
_{p}^{\gamma},\qquad(x,t)\in C^{-}.\label{segno1}%
\end{equation}
Assume that \eqref{1.17z} is false. Then for every $j\in\mathbf{N}$, there
exists a positive integer $k_{j}$ and functions $(u_{j},f_{j},g_{j}%
,\varphi_{j})\in{\mathcal{P}}_{0}(\mathcal{H},\alpha,C^{-},M_{1},M_{2},M_{3})$
such that $u_{j}(0,0)=\varphi_{j}(0,0)=0$ and
\begin{equation}
S_{k_{j}+1}^{-}(u_{j})>\max\left(  \frac{jM_{3}}{2^{(k_{j}+1)\gamma}}%
,\frac{S_{k_{j}}^{-}(u_{j})}{2^{\gamma}},\frac{S_{k_{j}-1}^{-}(u_{j})}{2^{2\gamma}%
},\dots,\frac{S_{0}^{-}(u_{j})}{2^{(k_{j}+1)\gamma}}\right)  ,\label{3.13-}%
\end{equation}
while (\ref{1.17z}) is true for $k<k_{j}.$ Using the definition of
$S_{k_{j}+1}^{-}$ in \eqref{1.17x} we see that there exists $(x_{j},t_{j})$ in
the closure of $C_{2^{-(k_{j}+1)}}^{-}$ such that $|u_{j}(x_{j},t_{j}%
)|=S_{k_{j}+1}^{-}(u_{j})$ for every $j\geq1$. Moreover from \eqref{segno1}
and (\ref{3.13-}) it follows that $u_{j}(x_{j},t_{j})>0$. Using \eqref{3.13-}
we can conclude, as $|u_{j}|\leq M_{1}$, that $j2^{-\gamma k_{j}}$ is bounded and
hence that $k_{j}\rightarrow\infty$ as $j\rightarrow\infty$.

\medskip\noindent\textbf{Step 2.} We define $(\tilde{x}_{j},\tilde{t}%
_{j})=\delta_{2^{k_{j}}}(x_{j},t_{j})$ and $\tilde{u}_{j}:C_{2^{k_{j}}}%
^{-}\longrightarrow\mathbf{R}$,
\[
\tilde{u}_{j}(x,t)=\frac{u_{j}(\delta_{2^{-k_{j}}}(x,t))}{S_{k_{j}+1}%
^{-}(u_{j})}.
\]
In particular, we note that $(\tilde{x}_{j},\tilde{t}_{j})$ belongs to the
closure of $C_{1/2}^{-}$ and
\begin{equation}
\tilde{u}_{j}(\tilde{x}_{j},\tilde{t}_{j})=1.\label{unotil}%
\end{equation}
Moreover we let $\tilde{\mathcal{H}}_{j}=\mathcal{H}_{2^{-k_{j}}}$, see
\eqref{3.4} for the exact definition of the scaled operator, and let
\begin{equation}
\tilde{f}_{j}(x,t)=2^{-2k_{j}}\frac{f_{j}(\delta_{2^{-k_{j}}}(x,t))}%
{S_{k_{j}+1}^{-}(u_{j})},\quad\tilde{g}_{j}(x,t)=\frac{g_{j}(\delta
_{2^{-k_{j}}}(x,t))}{S_{k_{j}+1}^{-}(u_{j})},\quad\tilde{\varphi}%
_{j}(x,t)=\frac{\varphi_{j}(\delta_{2^{-k_{j}}}(x,t))}{S_{k_{j}+1}^{-}(u_{j}%
)},\label{3.16-}%
\end{equation}
whenever $(x,t)\in C_{2^{k_{j}}}^{-}$. Then, using \eqref{riscal} we see that
$\tilde{u}_{j}$ solves
\[%
\begin{cases}
\max\{\tilde{\mathcal{H}}_{j}\tilde{u}_{j}-\tilde{f}_{j},\tilde{\varphi}%
_{j}-\tilde{u}_{j}\}=0, & \text{in}\ C_{2^{k_{j}}}^{-},\\
\tilde{u}_{j}=\tilde{g}_{j}, & \text{on}\ \partial_{p}C_{2^{k_{j}}}^{-}.
\end{cases}
\]
Moreover, for any $l\in\mathbf{N}$ we have that
\begin{equation}
\sup_{C_{2^{l}}^{-}}|\tilde{u}_{j}|=\frac{S_{k_{j}-l}^{-}(u_{j})}{S_{k_{j}%
+1}^{-}(u_{j})}\leq2^{(l+1)\gamma}\mbox{ whenever }k_{j}>l.\label{3.15-}%
\end{equation}
In particular, we can conclude that
\begin{equation}
(\tilde{u}_{j},\tilde{f}_{j},\tilde{u}_{j},\tilde{\varphi}_{j})\in
{\mathcal{P}}_{0}(\tilde{\mathcal{H}}_{j},\alpha,C_{2^{l}}^{-},\tilde{M}%
_{1}^{j},\tilde{M}_{2}^{j},\tilde{M}_{3}^{j}),\mbox{ and that $\tilde
u_j(0,0)=\tilde\varphi_j(0,0)$,}\label{3.17}%
\end{equation}
and, using \eqref{3.16-} and \eqref{3.15-}, we have that
\[
\tilde{M}_{1}^{j}\leq2^{(l+1)\gamma},\qquad\tilde{M}_{2}^{j}\leq2^{-2k_{j}}%
\frac{M_{2}}{S_{k_{j}+1}^{-}(u_{j})},\qquad\tilde{M}_{3}^{j}\leq
2^{\gamma(l-k_{j})}\frac{M_{3}}{S_{k_{j}+1}^{-}(u_{j})}.
\]
Now, by \eqref{3.13-} we see that in $C_{2^{l}}^{-}$
\begin{equation}
\lim_{j\rightarrow\infty}\tilde{M}_{2}^{j}=\lim_{j\rightarrow\infty}\tilde
{M}_{3}^{j}=0.\label{3.19-}%
\end{equation}

\medskip\noindent\textbf{Step 3.} In the following we let $l$ be a suitably
large positive integer to be specified later. We consider $j_{0}\in\mathbf{N}$
such that $k_{j}>2^{l}$ for $j\geq j_{0}$. We let $\hat{g}_{j}$ denote the
boundary values of $\tilde{u}_{j}$ on $\partial_{p}C_{2^{l}}^{-}$ and we let
$v_{j}$ and $\tilde{v}_{j}$ be such that
\[%
\begin{cases}
\tilde{\mathcal{H}}_{j}v_{j}=\Vert\tilde{f}_{j}\Vert_{L^{\infty}(C_{2^{l}}%
^{-})}\  & \text{in}\ C_{2^{l}}^{-},\\
v_{j}=\hat{g}_{j}\  & \text{on}\ \partial_{p}C_{2^{l}}^{-},
\end{cases}
\qquad\ \ \text{and \ \ \ \ }%
\begin{cases}
\tilde{\mathcal{H}}_{j}\tilde{v}_{j}=-\Vert\tilde{f}_{j}\Vert_{L^{\infty
}(C_{2^{l}}^{-})}\  & \text{in}\ C_{2^{l}}^{-},\\
\tilde{v}_{j}=\max\{\hat{g}_{j},\tilde{M}_{3}^{j}\}\  & \text{on}%
\ \partial_{p}C_{2^{l}}^{-}.
\end{cases}
\]
We shall prove that
\begin{equation}
v_{j}\leq\tilde{u}_{j}\leq\tilde{v}_{j}\mbox{ in }C_{2^{l}}^{-}.\label{3.21-}%
\end{equation}
The first inequality in \eqref{3.21-} follows from the comparison principle.
To prove the second one, we note that $\Vert\tilde{\varphi}_{j}\Vert
_{L^{\infty}}\leq\tilde{M}_{3}^{j}$ and then by the maximum principle
$\tilde{v}_{j}\geq\tilde{\varphi}_{j}$ in $C_{2^{l}}^{-}$. Furthermore,
\[
\tilde{\mathcal{H}}_{j}(\tilde{v}_{j}-\tilde{u}_{j})=-\Vert\tilde{f}_{j}%
\Vert_{L^{\infty}(C_{2^{l}}^{-})}+\tilde{f}_{j}\leq0\quad\text{in}\quad
D:=C_{2^{l}}^{-}\cap\{(x,t):\ \tilde{u}_{j}(x,t)>\tilde{\varphi}_{j}(x,t)\},
\]
and $\tilde{v}_{j}\geq\tilde{u}_{j}$ on $\partial_{p}D$. Hence, the second
inequality in \eqref{3.21-} follows from the maximum principle. Further, since
$\tilde{u}_{j}\geq\tilde{\varphi}_{j}$ by \eqref{3.17}, we can conclude that
$\hat{g}_{j}\geq-\tilde{M}_{3}^{j}$ in $C_{2^{l}}^{-}$. Now we use the maximum
principle to see that
\begin{equation}
\tilde{v}_{j}(x,t)-v_{j}(x,t)\leq\left(  \max\{0,\tilde{M}_{3}^{j}-\hat{g}%
_{j}\}+2\Vert\tilde{f}_{j}\Vert_{L^{\infty}(C_{2^{l}}^{-})}\right)
\leq2\left(  \tilde{M}_{3}^{j}+\tilde{M}_{2}^{j}\right) \label{3.27}%
\end{equation}
whenever $(x,t)\in C_{2^{l},1}^{-}$. We claim that there exists a positive
constant $c_{1}$ such that
\begin{equation}
\tilde{v}_{j}(x,t)\geq c_{1}\quad\mbox{for every }(x,t)\in C_{1/2}^{-},\ j\geq
j_{0}.\label{3.26}%
\end{equation}
Assuming \eqref{3.26} we use \eqref{3.27} and \eqref{3.21-} to conclude that
\[
\tilde{u}_{j}(0,0)\geq v_{j}(0,0)\geq\tilde{v}_{j}(0,0)-2\left(  \tilde{M}%
_{3}^{j}+\tilde{M}_{2}^{j}\right)  \geq c_{1}-2\left(  \tilde{M}_{3}%
^{j}+\tilde{M}_{2}^{j}\right)  ,
\]
and then, by \eqref{3.19-}, that $\tilde{u}_{j}(0,0)>0$ for $j$ suitably
large. This contradicts the assumption that $\tilde{u}_{j}(0,0)=\tilde
{\varphi}_{j}(0,0)=0$. Hence our original assumption was incorrect and the
proof of the lemma is complete. Based on the above it only remains to prove
\eqref{3.26}. Our proof of \eqref{3.26} is based on Lemma \ref{Harnack} and
Lemma \ref{t-2.19}. Below we will use the following short notation for $R_{1}
$, $R_{2}\geq1$,
\[
\partial_{p}^{+}C_{R_{1},R_{2}}^{-}=\partial_{p}C_{R_{1},R_{2}}^{-}%
\cap\{t>-R_{2}^{2}\},\ \text{\ \ \ \ }\partial_{p}^{-}C_{R_{1},R_{2}}%
^{-}=\partial_{p}C_{R_{1},R_{2}}^{-}\cap\{t=R_{2}^{2}\},
\]
where $C_{R_{1},R_{2}}^{-}=B_{d}(0,R_{1})\times(-R_{2}^{2},0]$. We write
$\tilde{v}_{j}=w_{j}+\tilde{w}_{j}+\hat{w}_{j}$ on $C_{2^{l},1}^{-}$ where
\[%
\begin{cases}
\tilde{\mathcal{H}}_{j}w_{j}=0\  & \text{in}\ C_{2^{l},1}^{-},\\
w_{j}=0\  & \text{on}\ \partial_{p}^{+}C_{2^{l},1}^{-},\\
w_{j}=\tilde{v}_{j}\  & \text{on}\ \partial_{p}^{-}C_{2^{l},1}^{-},
\end{cases}
\qquad%
\begin{cases}
\tilde{\mathcal{H}}_{j}\tilde{w}_{j}=0\  & \text{in}\ C_{2^{l},1}^{-},\\
\tilde{w}_{j}=\tilde{v}_{j}\  & \text{on}\ \partial_{p}^{+}C_{2^{l},1}^{-},\\
\tilde{w}_{j}=0\  & \text{on}\ \partial_{p}^{-}C_{2^{l},1}^{-},
\end{cases}
\text{ \ \ \ \ }%
\begin{cases}
\tilde{\mathcal{H}}_{j}\hat{w}_{j}=-\Vert\tilde{f}_{j}\Vert_{L^{\infty
}(C_{2^{l}}^{-})}\  & \text{in}\ C_{2^{l},1}^{-},\\
\hat{w}_{j}=0\  & \text{on}\ \partial_{p}C_{2^{l},1}^{-}.
\end{cases}
\]

By the maximum principle we see that
\begin{equation}
0\leq\hat{w}_{j}(x,t)\leq(t+1)\Vert\tilde{f}_{j}\Vert_{L^{\infty}(C_{2^{l}%
}^{-})}\leq\tilde{M}_{2}^{j},\label{ww}%
\end{equation}
whenever $(x,t)\in C_{2^{l},1}^{-}$. Hence, as $t\in\,(-1,0)$ we see that
$|\hat{w}_{j}(x,t)|\leq1/4$ in $C_{2^{l},1}^{-}$ if $j$ is sufficiently large.
Since
\begin{equation}\label{vj}
\Vert\tilde{v}_{j}\Vert_{L^{\infty}\left(  C_{2^{l}}^{-}\right)  }\leq %
\max\left\{  2^{(l+1)\gamma},\tilde{M}_{3}^{j}\right\}  +2^{2l}\tilde{M}_{2}^{j},
\end{equation}
by Lemma \ref{t-2.19} we find
\begin{equation}
\sup_{C_{1/2}^{-}}|\tilde{w}_{j}|\leq c2^{(Q+1)l}e^{-c2^{2l}}\sup
_{\partial_{p}^{+}C_{2^{l}}^{-}}|v|\leq c2^{(Q+1)l}e^{-c2^{2l}}\left(
\max\left\{  2^{(l+1)\gamma},\tilde{M}_{3}^{j}\right\}  +2^{2l}\tilde{M}_{2}%
^{j}\right)  .\label{foto}%
\end{equation}
In particular we note that the right hand side in this inequality tends to
zero as $l$ goes to infinity. Recalling that $\tilde{v}_{j}(\tilde{x}%
_{j},\tilde{t}_{j})\geq\tilde{u}_{j}(\tilde{x}_{j},\tilde{t}_{j})=1$, we can
conclude, by choosing $l$ suitably large, that
\[
w_{j}(\tilde{x}_{j},\tilde{t}_{j})\geq\frac{1}{2},\qquad j\geq j_{0}.
\]
Using this and the maximum principle we can conclude that there exists at
least one point $(\bar{x}_{j},\bar{t}_{j})\in\partial_{p}^{-}C_{2^{l},1}^{-}$
such that
\[
\tilde{v}_{j}(\bar{x}_{j},\bar{t}_{j})=w_{j}(\bar{x}_{j},\bar{t}_{j})\geq
\frac{1}{2},\qquad j\geq j_{0}.
\]
Thereafter we decompose $\tilde{v}_{j}=\check{v}_{j}+\hat{v}_{j}$ where
\[%
\begin{cases}
\tilde{\mathcal{H}}_{j}\check{v}_{j}=0\  & \text{in}\ C_{2^{l},2}^{-},\\
\check{v}_{j}=\tilde{v}_{j}\  & \text{on}\ \partial_{p}C_{2^{l},2}^{-},
\end{cases}
\qquad%
\begin{cases}
\tilde{\mathcal{H}}_{j}\hat{v}_{j}=-\Vert\tilde{f}_{j}\Vert_{L^{\infty
}(C_{2^{l}}^{-})}\  & \text{in}\ C_{2^{l},2}^{-},\\
\hat{v}_{j}=0\  & \text{on}\ \partial_{p}C_{2^{l},2}^{-}.
\end{cases}
\]
As in \eqref{ww} we see that $|\hat{v}_{j}(x,t)|\leq1/4$ in $C_{2^{l},2}^{-} $
if $j$ is sufficiently large and hence we can conclude that
\[
\check{v}_{j}(\bar{x}_{j},\bar{t}_{j})\geq\frac{1}{4},\qquad j\geq j_{0}.
\]
Using Lemma \ref{Harnack} we can therefore conclude that
\[
\inf_{C_{1/2}^{-}}\check{v}_{j}\geq\frac{1}{4\tilde{c}}.
\]
Further, since $\hat{v}_{j}\rightarrow0$ uniformly on $C_{2^{l},2}^{-}$ as $j$
goes to infinity, we conclude that
\[
\inf_{C_{1/2}^{-}}\tilde{v}_{j}\geq\inf_{C_{1/2}^{-}}\check{v}_{j}-\Vert
\hat{v}_{j}\Vert_{L^{\infty}}\geq\frac{1}{8\tilde{c}}%
\]
for any $j,$ suitably large. In particular, this proves \eqref{3.26} and hence
the proof of Lemma \ref{lem1} is complete in the case $m=0$.

\bigskip

Next, we consider the case when $\gamma=1+\alpha$\ and $m=1.$ In this
particular case we need to prove that (\ref{1.17z}) holds with $F=P_{1}%
^{(0,0)}\varphi$ and we assume that%
\[
(u,f,g,\varphi)\in\mathcal{P}_{1}(\mathcal{H},\alpha,C^{-},M_{1},M_{2},M_{3}).
\]
As in the previous case, we divide the argument in three steps, and the proof
will be similar to the one just carried out.

\textbf{Step 1. }By using that $\varphi(0,0)=0$ and that $\varphi
\in\mathcal{C}^{1,\alpha}(C^{-}),$ we note that%
\[
u(x,t)\geq\varphi(x,t)-\varphi(0,0)\geq F(x,t)-M_{3}||(x,t)||_{p}^{\gamma
},\text{ \ \ \ \ }(x,t)\in C^{-}.
\]
As previously, we assume that (\ref{1.17z}) is false. That is, for every
$j\in\mathbf{N}$ there exist a (smallest) positive integer $k_{j}$ and
$(u_{j},f_{j},g_{j},\varphi_{j})\in\mathcal{P}_{1}(\mathcal{H},\alpha
,C^{-},M_{1},M_{2},M_{3})$ such that $u_{j}(0,0)\geq0$ and%
\begin{equation}
S_{k_{j}+1}^{-}(u_{j}-F_{j})>\max\left(  \frac{jM_{3}}{2^{(k_{j}+1)\gamma}%
},\frac{S_{k_{j}}^{-}(u_{j}-F_{j})}{2^{\gamma}},\frac{S_{k_{j}-1}^{-}%
(u_{j}-F_{j})}{2^{2\gamma}},...,\frac{S_{0}^{-}(u_{j}-F_{j})}{2^{(k_{j}%
+1)\gamma}}\right)  .\label{dagis}%
\end{equation}
Here $F_{j}=P_{1}^{(0,0)}\varphi_{j}$ which in view of Theorem \ref{Taylor} is
given by%
\[
F_{j}(x,t)=\sum_{i=1}^{q}X_{i}\varphi_{j}(0,0)\cdot x_{i}.
\]
Using \eqref{1.17x} we see that there exists $(x_{j},t_{j})$ in the closure of
$C_{2^{-(k_{j}+1)}}^{-}$ such that $|u_{j}(x_{j},t_{j})-F_{j}(x_{j}%
,t_{j})|=S_{k_{j}+1}^{-}(u_{j})$ for every $j\geq1$. Moreover, by definition,
$u_{j}(x_{j},t_{j})-F_{j}(x_{j},t_{j})>0$. Since $u_{j}-F_{j}$ is bounded and
using (\ref{dagis}) we conclude that $j2^{-(k_{j}+1)\gamma}$ is bounded, and
hence, that $k_{j}\rightarrow\infty$ as $j\rightarrow\infty$.

\textbf{Step 2. }We define $(\tilde{x}_{j},\tilde{t}_{j})=\delta_{2^{k_{j}}%
}(x_{j},t_{j})$ and $\tilde{u}_{j}:C_{2^{k_{j}}}^{-}\longrightarrow\mathbf{R}%
$,
\[
\tilde{u}_{j}(x,t)=\frac{(u_{j}-F_{j})(\delta_{2^{-k_{j}}}(x,t))}{S_{k_{j}%
+1}^{-}(u_{j}-F_{j})}.
\]
As before, we note that $(\tilde{x}_{j},\tilde{t}_{j})$ belongs to the closure
of $C_{1/2}^{-}$ and that
\[
\tilde{u}_{j}(\tilde{x}_{j},\tilde{t}_{j})=1.
\]
Similarly we define
\[
\tilde{f}_{j}(x,t)=2^{-2k_{j}}\frac{(f_{j}-\mathcal{H}F_{j})(\delta
_{2^{-k_{j}}}(x,t))}{S_{k_{j}+1}^{-}(u_{j}-F_{j})},\quad\tilde{g}%
_{j}(x,t)=\frac{(g_{j}-F_{j})(\delta_{2^{-k_{j}}}(x,t))}{S_{k_{j}+1}^{-}%
(u_{j}-F_{j})},\quad\tilde{\varphi}_{j}(x,t)=\frac{(\varphi_{j}-F_{j}%
)\delta_{2^{-k_{j}}}(x,t))}{S_{k_{j}+1}^{-}(u_{j}-F_{j})}%
\]
and consider the scaled operator $\tilde{\mathcal{H}}_{j}=\mathcal{H}%
_{2^{-k_{j}}}$, see \eqref{3.4}. We also note that $\tilde{u}_{j}$ solves
\[
\left\{
\begin{array}
[c]{ll}%
\max\left\{  \tilde{\mathcal{H}}_{j}\tilde{u}_{j}-\tilde{f}_{j},\tilde
{\varphi}_{j}-\tilde{u}_{j}\right\}  =0 & \text{in }C_{2^{k_{j}}}^{-},\\
\tilde{u}_{j}=\tilde{g}_{j} & \text{on }\partial_{p}C_{2^{k_{j}}}^{-}.
\end{array}
\right.
\]
Moreover, for any $l\in\mathbf{N}$ with $l<k_{j}$ we have $(\tilde{u}%
_{j},\tilde{f}_{j},\tilde{u}_{j},\tilde{\varphi}_{j})\in{\mathcal{P}}%
_{1}(\tilde{\mathcal{H}}_{j},\alpha,C_{2^{l}}^{-},\tilde{M}_{1}^{j},\tilde
{M}_{2}^{j},\tilde{M}_{3}^{j})$. This time we prove our claim by using the
seminorm%
\begin{equation}
N_{3}^{j}:=\sum_{i=1}^{q}\left\vert X_{i}^{2^{-k_{j}}}\tilde{\varphi}%
_{j}(0,0)\right\vert +\underset{z,\zeta\in C_{2^{l}}^{-},\ z\neq\zeta}{\sup
}\frac{\left\vert \tilde{\varphi}_{j}(z)-\tilde{\varphi}_{j}(\zeta)-\sum
_{i=1}^{q}(z_{i}-\zeta_{i})X_{i}^{2^{-k_{j}}}\tilde{\varphi}_{j}%
(\zeta)\right\vert }{d_{p}(z,\zeta)^{\gamma}},\label{nlp}%
\end{equation}
where $X_{i}^{2^{-k_{j}}},\ i=1,2,...,q,$ are the scaled vector fields defined
as in Subsection \ref{SecFC}. By definition,
\[
N_{j}^{3}+||\tilde{\varphi}_{j}||_{\mathcal{C}^{0,\alpha}(C_{2^{l}}^{-})}%
\leq\tilde{M}_{3}^{j},
\]
where the H\"{o}lder norm is defined in terms of the scaled vector fields
$\{X_{i}^{2^{-k_{j}}}\}.$ A direct calculation shows that
\[
\mathcal{H}F_{j}(x,t)=\sum_{i=1}^{q}b_{i}(x,t)X_{i}(\sum_{j=1}^{q}X_{j}%
\varphi(0,0)\cdot x_{j})=\sum_{i,j=1}^{q}b_{i}(x,t)X_{j}\varphi(0,0)X_{j}%
x_{i}=\sum_{i,j=1}^{q}c_{ji}(x)b_{i}(x,t)X_{j}\varphi(0,0),
\]
and therefore
\[
\tilde{M}_{1}^{j}\leq2^{(l+1)\gamma},\ \ \ \tilde{M}_{2}^{j}\leq2^{-2k_{j}%
}\frac{M_{2}+c(\mathbf{c}_{\alpha},X)M_{3}}{S_{k_{j}+1}^{-}(u_{j}-F_{j}%
)}\text{ \ \ and \ \ }N_{3}^{j}\leq2^{(l-k_{j})\gamma}\frac{M_{3}}{S_{k_{j}%
+1}^{-}(u_{j}-F_{j})}.
\]
Note that the first two inequalities follow directly from the definition while
the third one follows from the fact that $X_{i}\tilde{\varphi}_{j}(0,0)=0$ and
the following calculation:%
\begin{align*}
& \left\vert \left(  \varphi_{j}-F_{j}\right)  (z)-\left(  \varphi_{j}%
-F_{j}\right)  (\zeta)-\sum_{i=1}^{q}(x_{i}-\xi_{i})X_{i}\left(  \varphi
_{j}-F_{j}\right)  (\zeta)\right\vert \\
& \leq||\varphi_{j}||_{\mathcal{C}^{1,\alpha}\left(  C_{2^{l}}^{-}\right)
}d_{p}(z,\zeta)^{1+\alpha}+\left\vert F_{j}(z)-F_{j}(\zeta)-\sum_{i=1}%
^{q}(x_{i}-\xi_{i})X_{i}F_{j}(\zeta)\right\vert \\
& =||\varphi_{j}||_{\mathcal{C}^{1,\alpha}\left(  C_{2^{l}}^{-}\right)  }%
d_{p}(z,\zeta)^{1+\alpha}.
\end{align*}
We also mention that the last equality follows from the definition of $F_{j}$
and the fact that $X_{i}x_{j}=\delta_{ij}$ when $X_{i}$ belongs to the first
layer of the Jacobian basis of $\mathbf{G}$. Therefore, by (\ref{dagis}),%
\[
\lim_{j\rightarrow\infty}\tilde{M}_{2}^{j}=\lim_{j\rightarrow\infty}N_{3}%
^{j}=0.
\]

\textbf{Step 3. }First of all, $||\tilde{\varphi}_{j}||_{L^{\infty}}\leq
N_{3}^{j}$, and this follows from the definition of $N_{3}^{j}$ in (\ref{nlp})
with $\zeta=0$. From here this step is completely analogous to Step 3 in the
first part of the proof, that is when $m=0$ and $\gamma=\alpha$, but with
$\tilde{M}_{3}^{j}$ replaced by $N_{3}^{j}$. This completes the proof when
$m=1,$ $\gamma=1+\alpha$.

\noindent It remains to consider the case when $m=2$\ and $\gamma=2$. In this
particular case we need to prove that (\ref{1.17z}) holds with $F=\varphi$ and
we assume that%
\[
(u,f,g,\varphi)\in\mathcal{P}_{2}(\mathcal{H},\alpha,C^{-},M_{1},M_{2},M_{3}).
\]
Below we will assume, without loss of generality, that $\varphi\equiv0.$
Indeed, since $\varphi\in\mathcal{C}^{2,\alpha}(C^{-}),$ $\mathcal{H}%
\varphi\in\mathcal{C}^{0,\alpha}(C^{-})$ and we can consider the modified
functions $\widehat{u}=u-\varphi$, $\widehat{f}=f-\mathcal{H}\varphi,$
$\widehat{g}=g-\varphi$ in $C^{-}.$ Note that
\[
(\widehat{u},\widehat{f},\widehat{g},0)\in\mathcal{P}_{2}(\mathcal{H}%
,\alpha,C^{-},\widehat{M}_{1},\widehat{M}_{2},0),
\]
where $\widehat{M}_{1}$ and $\widehat{M}_{2}$ (only) depend on $M_{1},\ M_{2}
$ and $M_{3}.\ $In this setting, we can now argue as before, but this time,
the contradictive assumption is that for every $j\in\mathbf{N}$ there exist a
positive integer $k_{j}$ and functions $(u_{j},f_{j},g_{j},0)\in
\mathcal{P}_{2}(\mathcal{H},\alpha,C^{-},M_{1},M_{2},0)$ such that%
\begin{equation}
S_{k_{j}+1}^{-}(u_{j})>\max\left(  \frac{j}{2^{(k_{j}+1)\gamma}}%
,\frac{S_{k_{j}}^{-}(u_{j})}{2^{\gamma}},\frac{S_{k_{j}-1}^{-}(u_{j}%
)}{2^{2\gamma}},...,\frac{S_{0}^{-}(u_{j})}{2^{(k_{j}+1)\gamma}}\right)
.\label{dagis2}%
\end{equation}
Using \eqref{1.17x} we see that there exists $(x_{j},t_{j})$ in the closure of
$C_{2^{-(k_{j}+1)}}^{-}$ such that $|u_{j}(x_{j},t_{j})|=S_{k_{j}+1}^{-}%
(u_{j})$ for every $j\geq1$. Moreover, by definition, $u_{j}(x_{j},t_{j})>0$.
Using (\ref{dagis2}) we can conclude, as $u_{j}$ is bounded, that
$j2^{-(k_{j}+1)\gamma}$ is bounded and hence that $k_{j}\rightarrow\infty$ as
$j\rightarrow\infty$. From now on, the proof follows along the same lines as
in the case $m=0$, $\gamma=\alpha$ and we omit the details. This completes the
proof.\hfill$\Box$

\subsection{Proof of Lemma \ref{lem}}

Assume $(u,f,g,\varphi)\in{\mathcal{\tilde{P}}}_{m}(\mathcal{H},\alpha
,C^{+},M_{1},M_{2},M_{3},M_{4})$. Then we first note that there exists a
constant $\eta_{m,\alpha}=\eta(\mathcal{H},m,\alpha,M_{1},M_{2},M_{3})$ such
that
\begin{equation}
\inf_{C_{r}^{+}}(u-g)\geq-\eta_{m,\alpha}r^{\gamma},\qquad r\in
(0,1),\label{low1}%
\end{equation}
where $\gamma=m+\alpha$ for $m=0,1$ and $\gamma=2$ for $m=2$. In fact, to
prove this we let $v$ be the solution to the Dirichlet problem
\eqref{e-obsobs} in the domain $D=C^{+}$ with right-hand side defined by $f$
and with boundary defined by $g$. Then
\[
(v,f,g)\in{\mathcal{D}}_{m}(\mathcal{H},\alpha,C^{+},M_{1},M_{2},M_{3}).
\]
Furthermore, since $v$ solves the Dirichlet problem while $u$ solves the
obstacle problem with the same boundary data we have $u\geq v$ on the closure
of $C^{+}$ and hence \eqref{low1} follows directly if we apply Lemma
\ref{lemD} to $v$. We next note, as in the proof of Lemma \ref{lemD}, that it
suffices to prove that
\begin{equation}
\sup_{C_{r}^{+}}\left\vert u-P_{m}^{(0,0)}g\right\vert \leq cr^{\gamma},\qquad
r\in(0,1),\label{a1}%
\end{equation}
where $\gamma=m+\alpha$ for $m=0,1$ and
\begin{equation}
\sup_{C_{r}^{+}}\left\vert u-g\right\vert \leq cr^{\gamma},\qquad r\in
(0,1),\label{b1}%
\end{equation}
where $\gamma=2$ for $m=2$. Furthermore, as in the proof of Lemma \ref{lemD},
we can without loss of generality assume that $P_{m}^{(0,0)}g\equiv0$ for
$m=0,1$ and that $g\equiv0$ for $m=2$. To prove \eqref{a1} and \eqref{b1}
after these preliminary problem reduction steps we define $S_{k}^{+}(u)$ as in
\eqref{1.17x+} and we intend to prove that there exists a positive constant
$\tilde{c}=\tilde{c}\left(  \mathcal{H},m,\alpha,M_{1},M_{2},M_{3}%
,M_{4}\right)  $ such that \eqref{1.17z+} holds for all $k\in\mathbf{N}$ with
$F$ as defined in \eqref{lex2} (which by our problem reduction reads
$F\equiv0$). Again, if \eqref{1.17z+} holds then
\[
S_{k}^{+}(u-F)\leq\frac{\tilde{c}}{2^{k\gamma}}%
\]
and hence \eqref{a1} and \eqref{b1}, and Lemma \ref{lem}, follow.

We first consider the case $m=0$\ and prove \eqref{1.17z+} with $\gamma=\alpha
$. In particular, we assume that
\[
(u,f,g,\varphi)\in{\mathcal{\tilde{P}}}_{0}(\mathcal{H},\alpha,C^{+}%
,M_{1},M_{2},M_{3},M_{4}),
\]
and, as in the proof of Lemma \ref{lem1} and following the lines of
\cite{NPP}, we divide the argument into three steps.

\medskip\noindent\textbf{Step 1. } We first note, using \eqref{low1}, that
\begin{equation}
u(x,t)\geq-\left(  \eta_{m,\alpha}+M_{3}\right)  \Vert(x,t)\Vert_{p}^{\gamma
},\qquad(x,t)\in C^{+}.\label{segno}%
\end{equation}
Assume that \eqref{1.17z+} is false. Then for every $j\in\mathbf{N}$, there
exists a positive integer $k_{j}$ and functions $(u_{j},f_{j},g_{j}%
,\varphi_{j})\in{\mathcal{\tilde{P}}}_{0}(\mathcal{H},\alpha,C^{+},M_{1}%
,M_{2},M_{3},M_{4})$ such that $u_{j}(0,0)=0\geq\varphi_{j}(0,0)$ and
\begin{equation}
S_{k_{j}+1}^{+}(u_{j})>\max\left(  \frac{j\left(  \eta_{m,\alpha}%
+M_{3}\right)  }{2^{(k_{j}+1)\gamma}},\frac{S_{k_{j}}^{+}(u_{j})}{2^{\gamma}}%
,\frac{S_{k_{j}-1}^{+}(u_{j})}{2^{2\gamma}},\dots,\frac{S_{0}^{+}(u_{j}%
)}{2^{(k_{j}+1)\gamma}}\right)  .\label{3.13}%
\end{equation}
Using the definition of $S_{k_{j}+1}^{+}$ in \eqref{1.17x+} we see that there
exists $(x_{j},t_{j})$ in the closure of $C_{2^{-(k_{j}+1)}}^{+}$ such that
$|u_{j}(x_{j},t_{j})|=S_{k_{j}+1}^{+}(u_{j})$ for every $j\geq1$. Moreover
from \eqref{segno} and (\ref{3.13}) it follows that $u_{j}(x_{j},t_{j})>0$.
Using \eqref{3.13} we can conclude, as $|u_{j}|\leq M_{1}$, that
$j2^{-\gamma k_{j}}$ is bounded and hence that $k_{j}\rightarrow\infty$ as
$j\rightarrow\infty$.

\medskip\noindent\textbf{Step 2.} We define $(\tilde{x}_{j},\tilde{t}%
_{j})=\delta_{2^{k_{j}}}((x_{j},t_{j}))$ and $\tilde{u}_{j}:C_{2^{k_{j}}}%
^{+}\longrightarrow\mathbf{R}$,
\[
\tilde{u}_{j}(x,t)=\frac{u_{j}(\delta_{2^{-k_{j}}}(x,t))}{S_{k_{j}+1}%
^{+}(u_{j})}.
\]
Note that $(\tilde{x}_{j},\tilde{t}_{j})$ belongs to the closure of
$C_{1/2}^{+}$ and
\begin{equation}
\tilde{u}_{j}(\tilde{x}_{j},\tilde{t}_{j})=1.\label{unotilp}%
\end{equation}
Moreover, we let $\tilde{\mathcal{H}}_{j}=\mathcal{H}_{2^{-k_{j}}}$, see
\eqref{3.4} for the exact definition of this scaled operator, and
\begin{equation}
\tilde{f}_{j}(x,t)=2^{-2k_{j}}\frac{f_{j}(\delta_{2^{-k_{j}}}(x,t))}%
{S_{k_{j}+1}^{+}(u_{j})},\quad\tilde{g}_{j}(x,t)=\frac{g_{j}(\delta
_{2^{-k_{j}}}(x,t))}{S_{k_{j}+1}^{+}(u_{j})},\quad\tilde{\varphi}%
_{j}(x,t)=\frac{\varphi_{j}(\delta_{2^{-k_{j}}}(x,t))}{S_{k_{j}+1}^{+}(u_{j}%
)}\label{3.16}%
\end{equation}
whenever $(x,t)\in C_{2^{k_{j}}}^{+}$. Then, using \eqref{riscal}, we see
that
\[%
\begin{cases}
\max\{\tilde{\mathcal{H}}_{j}\tilde{u}_{j}-\tilde{f}_{j},\tilde{\varphi}%
_{j}-\tilde{u}_{j}\}=0, & \text{in}\ C_{2^{k_{j}}}^{+},\\
\tilde{u}_{j}=\tilde{g}_{j}, & \text{on}\ \partial_{p}C_{2^{k_{j}}}^{+}.
\end{cases}
\]
In the following we let $l\in\mathbf{N}$ be fixed and to be specified below.
From \eqref{3.13} it follows that
\begin{equation}
\sup_{C_{2^{l}}^{+}}|\tilde{u}_{j}|=\frac{S_{k_{j}-l}^{+}(u_{j})}{S_{k_{j}%
+1}^{+}(u_{j})}\leq2^{(l+1)\gamma}\ \mbox{ whenever }\ k_{j}>l,\label{3.15}%
\end{equation}
and that
\[
(\tilde{u}_{j},\tilde{f}_{j},\tilde{u}_{j},\tilde{\varphi}_{j})\in
\widetilde{\mathcal{P}}_{0}(\tilde{\mathcal{H}}_{j},\alpha,C_{2^{l}}%
^{+},\tilde{M}_{1}^{j},\tilde{M}_{2}^{j},\tilde{M}_{3}^{j},\tilde{M}_{4}^{j}),
\]
for some $\tilde{M}_{1}^{j},\tilde{M}_{2}^{j},\tilde{M}_{3}^{j},\tilde{M}%
_{4}^{j}$. Furthermore, using \eqref{3.16}, \eqref{3.15} and the continuity of
$\tilde{f}_{j}$, we have
\begin{equation}
\tilde{M}_{1}^{j}\leq2^{(l+1)\gamma},\qquad\tilde{M}_{2}^{j}\leq2^{-2k_{j}}%
\frac{M_{2}}{S_{k_{j}+1}^{+}(u_{j})}.\label{3.18p}%
\end{equation}
Moreover, we let
\[
m_{j}=\max\bigg\{\left\Vert \tilde{g}_{j}\right\Vert _{L^{\infty}(C_{2^{l}%
}^{+})},\sup_{C_{2^{l}}^{+}}\tilde{\varphi}_{j}\bigg\}.
\]
Then, using \eqref{3.13} and the $\mathcal{C}_{X}^{0,\alpha}$-regularity of
$g_{j}$ and $\varphi_{j}$, we see that
\begin{equation}
\lim_{j\rightarrow\infty}\tilde{M}_{2}^{j}=\lim_{j\rightarrow\infty}%
m_{j}=0.\label{3.19}%
\end{equation}
Note that we can not ensure the decay of $\tilde{M}_{4}^{j}$, as
$j\rightarrow\infty$, as we only know that $\tilde{\varphi}_{j}(0,0)\leq0$.

\medskip\noindent\textbf{Step 3.} Below we will choose $l$ suitably large to
find a contradiction and we will use the following short notation for $R$,
$R_{1}$, $R_{2}\geq1$,
\begin{align}
\partial_{p}^{+}C_{{R}}^{+}  & =\partial_{p}C_{{R}}^{+}\cap\{t>0\},\ \partial
_{p}^{-}C_{{R}}^{+}=\partial_{p}C_{{R}}^{+}\cap\{t=0\},\nonumber\\
\partial_{p}^{+}C_{R_{1},R_{2}}^{+}  & =\partial_{p}C_{R_{1},R_{2}}^{+}%
\cap\{t>0\},\ \partial_{p}^{-}C_{R_{1},R_{2}}^{+}=\partial_{p}C_{R_{1},R_{2}%
}^{+}\cap\{t=0\}.\label{bs}%
\end{align}
Above $C_{R_{1},R_{2}}^{+}$ denotes $B_{d}(0,R_{1})\times(0,R_{2}^{2}]$. We
now consider $j_{0}\in\mathbf{N}$ such that $k_{j}>2^{l}$ for $j\geq j_{0}$
and we let $\tilde{v}_{j}$ be the solution to
\[%
\begin{cases}
\tilde{\mathcal{H}}_{j}\tilde{v}_{j}=-\Vert\tilde{f}_{j}\Vert_{L^{\infty
}(C_{2^{l}}^{+})}\  & \text{in}\ C_{2^{l}}^{+},\\
\tilde{v}_{j}=\tilde{M}_{1}^{j}\  & \text{on}\ \partial_{p}^{+}C_{2^{l}}%
^{+},\\
\tilde{v}_{j}=m_{j}\  & \text{on}\ \partial_{p}^{-}C_{2^{l}}^{+}.
\end{cases}
\]
To continue we prove that
\begin{equation}
\tilde{u}_{j}\leq\tilde{v}_{j}\mbox{ in }C_{2^{l}}^{+},\label{3.21}%
\end{equation}
and that this contradicts (\ref{unotilp}). Indeed, by the maximum principle we
have $\tilde{v}_{j}\geq m_{j}\geq\tilde{\varphi}_{j}$ in $C_{2^{l}}^{+}$.
Furthermore
\[
\tilde{\mathcal{H}}_{j}(\tilde{v}_{j}-\tilde{u}_{j})=-\Vert\tilde{f}_{j}%
\Vert_{L^{\infty}(C_{2^{l}}^{+})}+\tilde{f}_{j}\leq0\quad\text{in}\quad
D:=C_{2^{l}}^{+}\cap\{(x,t):\ \tilde{u}_{j}(x,t)>\tilde{\varphi}_{j}(x,t)\},
\]
and $\tilde{v}_{j}\geq\tilde{u}_{j}$ on $\partial_{p}D$. Hence \eqref{3.21}
follows from the maximum principle. Next, we show that \eqref{3.21}
contradicts \eqref{unotilp}. We write $\tilde{v}_{j}=w_{j}+\tilde{w}_{j}%
+\hat{w}_{j}$ on $C_{2^{l},1}^{+}$ where
\[%
\begin{cases}
\tilde{\mathcal{H}}_{j}w_{j}=0\  & \text{in}\ C_{2^{l},1}^{+},\\
w_{j}=0\  & \text{on}\ \partial_{p}^{+}C_{2^{l},1}^{+},\\
w_{j}=\tilde{v}_{j}\  & \text{on}\ \partial_{p}^{-}C_{2^{l},1}^{+},
\end{cases}
\qquad%
\begin{cases}
\tilde{\mathcal{H}}_{j}\tilde{w}_{j}=0\  & \text{in}\ C_{2^{l},1}^{+},\\
\tilde{w}_{j}=\tilde{v}_{j}\  & \text{on}\ \partial_{p}^{+}C_{2^{l},1}^{+},\\
\tilde{w}_{j}=0\  & \text{on}\ \partial_{p}^{-}C_{2^{l},1}^{+},
\end{cases}
\text{ \ \ \ \ }%
\begin{cases}
\tilde{\mathcal{H}}_{j}\hat{w}_{j}=-\Vert\tilde{f}_{j}\Vert_{L^{\infty
}(C_{2^{l}}^{-})}\  & \text{in}\ C_{2^{l},1}^{+},\\
\hat{w}_{j}=0\  & \text{on}\ \partial_{p}C_{2^{l},1}^{+}.
\end{cases}
\]
By the maximum principle we first see that
\begin{equation}
w_{j}\leq m_{j}\quad\mbox{ in }\ C_{2^{l},1}^{+},\label{w1}%
\end{equation}
and that
\begin{equation}
\Vert\hat{w}_{j}\Vert_{L^{\infty}\left(  C_{2^{l},1}^{+}\right)  }\leq
\Vert\tilde{f}_{j}\Vert_{L^{\infty}(C_{2^{l}}^{+})}\leq\tilde{M}_{2}%
^{j}.\label{w2}%
\end{equation}

\noindent We next use Lemma \ref{t-2.19} in the cylinder $C_{2^{l},1}^{+}$,
and by \eqref{3.18p} we see that
\begin{equation}
\sup_{C^{+}}\tilde{w}_{j}\leq c2^{(Q+1)l}e^{-c4^{l}}\sup_{\partial_{p}%
^{+}C_{2^{l},1}^{+}}\tilde{v}_{j}\leq c2^{(Q+1)l}e^{-c4^{l}}\tilde{M}_{1}%
^{j}\leq c2^{(Q+1)l}e^{-c4^{l}}2^{(l+1)\gamma},\label{w3}%
\end{equation}
and, in particular, we note that the right hand side in this inequality can be
made arbitrarily small by choosing $l$ large enough, independent of $j$.
Combining \eqref{w1}, \eqref{w2} and \eqref{w3} we conclude that, for a
suitably large $l$ and $j_{0}$, we have
\[
\sup_{C^{+}}\tilde{v}_{j}\leq\frac{1}{2}\qquad\text{for any }j\geq j_{0},
\]
which contradicts (\ref{unotilp}) and \eqref{3.21}. This completes the proof
in the case $m=0$.

It remains to consider the cases $\gamma=m+\alpha$\ for $m=1$\ and $\gamma
=2$\ for $m=2$, which can be proved in complete analogy with step 1-3 above
only with a slight change of motivation in (\ref{3.19}). We omit the
details.\hfill$\Box$

\subsection{Proof of Lemma \ref{lem1+}}

Using Lemma \ref{lem1} we only need to prove the statement on $C_{r}^{+}.$ We
also note, as in the proof of Lemma \ref{lemD}, that using the triangle
inequality and Lemma \ref{polyy++} it suffices to prove that
\begin{equation}
\sup_{C_{r}^{+}(0,0)}\left\vert u-P_{m}^{(0,0)}\varphi\right\vert \leq
cr^{\gamma},\qquad r\in(0,1),\label{pisa1}%
\end{equation}
for $\gamma=m+\alpha$ for $m=0,1$ and
\begin{equation}
\sup_{C_{r}^{+}(0,0)}\left\vert u-\varphi\right\vert \leq cr^{\gamma},\qquad
r\in(0,1),\label{pisa2}%
\end{equation}
for $\gamma=2$ for $m=2.$ Furthermore, as in the proof of Lemma \ref{lemD}, we
can without loss of generality assume that $P_{m}^{(0,0)}\varphi\equiv0$ for
$m=0,1$ and that $\varphi\equiv0$ for $m=2$. To prove (\ref{pisa1}) and
(\ref{pisa2}) after these preliminary problem reduction steps we define
$S_{k}^{+}(u)$ as in \eqref{1.17x+} and we intend to prove that there exists a
positive constant $\tilde{c}=\tilde{c}\left(  \mathcal{H},m,\alpha,M_{1}%
,M_{2},M_{3},M_{4}\right)  $ such that \eqref{1.17z+} holds for all
$k\in\mathbf{N}$ with $F$ as defined in \eqref{lex2}, which again, by our
problem reduction reads $F\equiv0$. If \eqref{1.17z+} holds then
\begin{equation}
S_{k}^{+}(u-F)\leq\frac{\tilde{c}}{2^{k\gamma}}\label{pisafalse}%
\end{equation}
and hence (\ref{pisa1}), (\ref{pisa2}) and Lemma \ref{lem1+} follows.

We first consider the case $m=0$\ and prove \eqref{1.17z+} with $\gamma
=\alpha$\emph{.} In particular, we assume that
\[
(u,f,g,\varphi)\in\mathcal{P}_{0}(\mathcal{H},\alpha,C^{+},M_{1},M_{2},M_{3}),
\]
and we will proceed in a way similar to the proof Lemma \ref{lem}. That is, we
assume that (\ref{pisafalse}) is false. Then for every $j\in\mathbf{N}$ there
exist a positive integer $k_{j}$ and $(u_{j},f_{j},g_{j},\varphi_{j}%
)\in\mathcal{P}_{0}(\mathcal{H},\alpha,C^{+},M_{1},M_{2},M_{3})$ such that
$u_{j}(0,0)=\varphi_{j}(0,0)=0$ and%
\[
S_{k_{j}+1}^{+}(u_{j})>\max\left(  \frac{jM_{3}}{2^{(k_{j}+1)\gamma}}%
,\frac{S_{k_{j}}^{+}(u_{j})}{2^{\gamma}},\frac{S_{k_{j}-1}^{+}(u_{j}%
)}{2^{2\gamma}},...,\frac{S_{0}^{+}(u_{j})}{2^{(k_{j}+1)\gamma}}\right)  .
\]
Furthermore, there exists $(x_{j},t_{j})$ in the closure of $C_{2^{-(k_{j}%
+1)}}^{+}$ such that $|u_{j}(x_{j},t_{j})|=S_{k_{j}+1}^{+}(u_{j})$ for every
$j\geq1.$ As in Step 2 in the proof of Lemma \ref{lem}, we define $(\tilde
{x}_{j},\tilde{t}_{j})=\delta_{2^{k_{j}}}((x_{j},t_{j}))$ and $\tilde{u}%
_{j}:C_{2^{k_{j}}}^{+}\rightarrow\mathbf{R}$,
\[
\tilde{u}_{j}(x,t)=\frac{u_{j}(\delta_{2^{-k_{j}}}(x,t))}{S_{k_{j}+1}%
^{+}(u_{j})}.
\]
Note that $(\tilde{x}_{j},\tilde{t}_{j})$ belongs to the closure of
$C_{1/2}^{+}$ and
\[
\tilde{u}_{j}(\tilde{x}_{j},\tilde{t}_{j})=1.
\]
Also in this case we can define functions $\tilde{f}_{j},\tilde{\varphi}_{j}$
similarly to $\tilde{u}_{j},$ see (\ref{3.16}), and then
\[
(\tilde{u}_{j},\tilde{f}_{j},\tilde{u}_{j},\tilde{\varphi}_{j})\in
\mathcal{P}_{0}(\tilde{\mathcal{H}}_{j},\alpha,C_{2^{l}}^{+},\tilde{M}_{1}%
^{j},\tilde{M}_{2}^{j},\tilde{M}_{3}^{j}),
\]
where%
\[
\tilde{M}_{1}^{j}\leq2^{(l+1)\gamma},\qquad\tilde{M}_{2}^{j}\leq2^{-2k_{j}%
}\frac{M_{2}}{S_{k_{j}+1}^{-}(u_{j})},\qquad\tilde{M}_{3}^{j}\leq
2^{(l-k_{j})\gamma}\frac{M_{3}}{S_{k_{j}+1}^{-}(u_{j})}.
\]
We now intend to complete the argument by a contradiction. We fix a suitable
positive integer $l,$ to be specified below, and we consider $j_{0}%
\in\mathbf{N}$ such that $k_{j}>2^{l}$ for $j>j_{0}.$ Then\ as in Lemma
\ref{lem1} we prove that%
\begin{equation}
v_{j}\leq\tilde{u}_{j}\leq\tilde{v}_{j},\text{ \ \ in }C_{2^{l}}%
^{+},\label{sara1}%
\end{equation}
where $v_{j}$ and $\tilde{v}_{j}$ are defined by
\begin{equation}%
\begin{cases}
\tilde{\mathcal{H}}_{j}v_{j}=\Vert\tilde{f}_{j}\Vert_{L^{\infty}(C_{2^{l}}%
^{+})}\  & \text{in}\ C_{2^{l}}^{+},\\
v_{j}=\hat{g}_{j}\  & \text{on}\ \partial_{p}C_{2^{l}}^{+},
\end{cases}
\qquad\text{respectively \ \ \ }%
\begin{cases}
\tilde{\mathcal{H}}_{j}\tilde{v}_{j}=-\Vert\tilde{f}_{j}\Vert_{L^{\infty
}(C_{2^{l}}^{+})}\  & \text{in}\ C_{2^{l}}^{+},\\
\tilde{v}_{j}=\max\{\hat{g}_{j},\tilde{M}_{3}^{j}\}\  & \text{on}%
\ \partial_{p}C_{2^{l}}^{+},
\end{cases}
\label{defv}%
\end{equation}
for some function $\hat{g}_{j}\ $which coincides with $\tilde{u}_{j}$ on
$\partial_{p}C_{2^{l}}^{+}$. Moreover, we note that
\begin{equation}
\lim_{j\rightarrow\infty}\tilde{M}_{3}^{j}=0=\lim_{j\rightarrow\infty}%
||\hat{g}_{j}||_{L^{\infty}(\partial_{p}C_{2^{l}}^{+}\cap\{t=0\})}%
\ =0.\label{defbound}%
\end{equation}
We next show that (\ref{sara1})-(\ref{defbound}) leads to a contradiction. To
do this we recall the notation introduced in (\ref{bs}) and write $\tilde
{v}_{j}=w_{j}+\tilde{w}_{j}+\hat{w}_{j}$ on $C_{2^{l},1}^{+}$ where
\[%
\begin{cases}
\tilde{\mathcal{H}}_{j}w_{j}=0\  & \text{in}\ C_{2^{l},1}^{+},\\
w_{j}=0\  & \text{on}\ \partial_{p}^{+}C_{2^{l},1}^{+},\\
w_{j}=\tilde{v}_{j}\  & \text{on}\ \partial_{p}^{-}C_{2^{l},1}^{+},
\end{cases}
\qquad%
\begin{cases}
\tilde{\mathcal{H}}_{j}\tilde{w}_{j}=0\  & \text{in}\ C_{2^{l},1}^{+},\\
\tilde{w}_{j}=\tilde{v}_{j}\  & \text{on}\ \partial_{p}^{+}C_{2^{l},1}^{+},\\
\tilde{w}_{j}=0\  & \text{on}\ \partial_{p}^{-}C_{2^{l},1}^{+},
\end{cases}
\text{ \ \ \ \ }%
\begin{cases}
\tilde{\mathcal{H}}_{j}\hat{w}_{j}=-\Vert\tilde{f}_{j}\Vert_{L^{\infty
}(C_{2^{l}}^{-})}\  & \text{in}\ C_{2^{l},1}^{+},\\
\hat{w}_{j}=0\  & \text{on}\ \partial_{p}C_{2^{l},1}^{+}.
\end{cases}
\]
Arguing as above (\ref{foto}),%
\[
||\hat{w}_{j}||_{L^{\infty}(C_{1/2}^{+})}\leq\tilde{M}_{2}^{j},\ \ \ ||\tilde
{w}_{j}||_{L^{\infty}(C_{1/2}^{+})}\leq c2^{(Q+1)l}e^{-c4^{l}}\left(
\max\{2^{(l+1)\gamma},\tilde{M}_{3}^{j}\}+2^{2l}\tilde{M}_{2}^{j}\right)  .
\]
Since $\tilde{v}_{j}(\tilde{x}_{j},\tilde{t}_{j})\geq\tilde{u}_{j}(\tilde
{x}_{j},\tilde{t}_{j})=1$ we can conclude, by choosing $l$ large enough, that
both $||\hat{w}_{j}||_{L^{\infty}(C_{1/2}^{+})}$ and $||\tilde{w}%
_{j}||_{L^{\infty}(C_{1/2}^{+})}$ tend to zero as $j\rightarrow\infty$, so
that%
\[
w_{j}(\tilde{x}_{j},\tilde{t}_{j})\geq\frac{1}{2}\text{ \ \ for }j\geq j_{0}.
\]
Now, applying the maximum principle we can conclude that there exists at least
one point $(\overline{x},\overline{t})\in\partial_{p}^{-}C_{2^{l},1}^{+}$ such
that%
\[
\tilde{v}_{j}(\overline{x},\overline{t})=w_{j}(\overline{x},\overline{t}%
)\geq\frac{1}{2}\text{ \ \ for }j\geq j_{0},
\]
which contradicts (\ref{defv}) and (\ref{defbound}). Hence (\ref{pisa1}) holds
for $m=0,$ $\gamma=\alpha$.

To prove Lemma\emph{\ \ref{lem1+} }when\emph{\ }$m=1,\ \gamma=1+\alpha
$\emph{\ }or when\emph{\ }$m=2,\ \gamma=2$ we can use the same arguments as
above and we omit the details.\hfill$\Box$

\section{Proof of Theorem \ref{th-int}}

The proof might be somewhat difficult to get an overview of although
straightforward, since we have to divide it into different cases along the
way. Therefore we will deal with different choices of $m$, i.e., $m=0,1,2$, in
the subsections below.

\subsection{Proof of Theorem \ref{th-int} when $m=0,$ $\gamma=\alpha$}

We will divide the proof into different cases stemming from geometric
assumptions. To explain this further we take $R\in(0,1)$ and a constant
$C_{1}>1$ such that $C_{2R}\subseteq C_{2C_{1}R}(\hat{x},\hat{t})\subseteq C$
for all $(\hat{x},\hat{t})\in C_{2R}$ and we define $\mathcal{F}=\overline
{C}_{2R}\cap\{(x,t):u(x,t)=\varphi(x,t)\}.$ Our goal is to prove that there
exists a constant $\widehat{c}=\widehat{c}(\mathcal{H},\alpha,f,\varphi)$ such
that%
\begin{equation}
\underset{(x,t)\neq(\hat{x},\hat{t})}{\sup_{(x,t),(\hat{x},\hat{t})\in C_{R}}%
}\frac{|u(x,t)-u(\hat{x},\hat{t})|}{d_{p}((x,t),(\hat{x},\hat{t}))^{\alpha}%
}\leq\widehat{c}.\label{vinter}%
\end{equation}
Before we continue, we note that if $\mathcal{F=\emptyset}$ then $u$ is a
solution to the Dirichlet problem and we can apply Lemma \ref{t-schauder}
and\ the result follows immediately. On the other hand, if $(\hat{x},\hat
{t})\in C_{2R}\cap\mathcal{F}$, then we can translate both the operator and
the functions so that $(\hat{x},\hat{t})$ is our new origin and after this
reformulation we can use Lemma \ref{lem1+} to obtain (\ref{vinter}), but this
time for $(x,t)\in C_{2C_{1}R}(\hat{x},\hat{t}).$ Obviously the analogue of
this estimate holds whenever $(x,t)\in C_{2R}\cap\mathcal{F}$. In view of this
it remains to consider the case when $(x,t),(\hat{x},\hat{t})\in
C_{R}\backslash\mathcal{F}.$ When we continue we divide the proof in two cases
depending on how far apart the points $(x,t),(\hat{x},\hat{t})$ are compared
with the distance from $(x,t)$ to $\mathcal{F}$. For this reason we define
$r=d_{p}((x,t),\mathcal{F}):=\inf\{d_{p}((x,t),(\xi,\tau)):(\xi,\tau
)\in\mathcal{F}\}$ and for fixed $(x,t)$ we let $(\tilde{x},\tilde{t}%
)\in\mathcal{F}$ be such a minimizing point. 

\textbf{Case 1:} Assume that $(\hat{x},\hat{t})\in C_{R}\backslash
C_{r/2}(x,t)$, then $d_{p}((x,t),(\hat{x},\hat{t}))$ $>c_{0}r$ for some
constant $c_{0}.$ Using the triangle inequality we get%
\begin{equation}
d_{p}((\hat{x},\hat{t}),(\tilde{x},\tilde{t}))\leq\left(  1+\frac{1}{c_{0}%
}\right)  d_{p}((x,t),(\hat{x},\hat{t}))\label{t-triang}%
\end{equation}
and by using (\ref{vinter}) we get
\begin{align*}
|u(x,t)-u(\hat{x},\hat{t})|  & \leq|u(x,t)-u(\tilde{x},\tilde{t}%
)|+|u(\tilde{x},\tilde{t})-u(\hat{x},\hat{t})|\\
& \leq\widehat{c}\left(  d_{p}((x,t),(\hat{x},\hat{t}))^{\alpha}+d_{p}%
((\hat{x},\hat{t}),(\tilde{x},\tilde{t}))^{\alpha}\right) \\
& \leq c_{1}d_{p}((x,t),(\hat{x},\hat{t}))^{\alpha}.
\end{align*}
As we are about to prove (\ref{vinter}) we remark that we may use
(\ref{vinter}) in the computation above since $(\tilde{x},\tilde{t})\in
C_{2R}\cap\mathcal{F}$, and in that particular case we have proved that
(\ref{vinter}) is valid. Hence, (\ref{vinter}) holds also in this case.

\textbf{Case 2: }Assume that $(\hat{x},\hat{t})\in C_{r/2}(x,t).$ Then we
define a new function%
\[
v(\hat{x},\hat{t}):=u(\hat{x},\hat{t})-u(\tilde{x},\tilde{t}),
\]
and we remark that $u(\tilde{x},\tilde{t})$ is to be treated as a constant.
Then, since $(\tilde{x},\tilde{t})\in C_{2R}\cap\mathcal{F}$ we can use
(\ref{vinter}) once more and we see that%
\[
|v(\hat{x},\hat{t})|\leq cd_{p}((\hat{x},\hat{t}),(\tilde{x},\tilde
{t}))^{\alpha}\leq c_{2}r^{\alpha}.
\]
To proceed, we define the function%
\[
w(y,s):=\frac{v^{r,(x,t)}(y,s)}{r^{\alpha}}%
\]
where $v^{r,(x,t)}$ is defined as in (\ref{and-e1}). For $(y,s)\in C_{1/2}$
the following holds:

\begin{description}
\item[$\mathit{i)}$] $|w(y,s)|\leq c_{3}$ for some constant $c_{3},$

\item[$\mathit{ii)}$] $\mathcal{H}^{r,(x,t)}w(y,s)=r^{2-\alpha}f^{r,(x,t)}%
(y,s).$
\end{description}

\noindent The Schauder estimate in Lemma \ref{t-schauder} then implies that%
\[
|w(y,s)-w(0,0)|\leq c_{4}d_{p}((0,0),(y,s))^{\alpha}%
\]
for some positive constant $c_{4}$. By the definition of $w$ this concludes
the proof in the case $m=0$.

\subsection{\bigskip Proof of Theorem \ref{th-int} when $m=1, $ $\gamma
=1+\alpha$}

As before we take $R\in(0,1)$ and $C_{1}>1$ such that $C_{2R}\subseteq
C_{2C_{1}R}(\hat{x},\hat{t})\subseteq C$ for all $(\hat{x},\hat{t})\in C_{2R}
$ and we define $\mathcal{F}=\overline{C}_{2R}\cap\{(x,t):u(x,t)=\varphi
(x,t)\}.$ We assume that $\mathcal{F\neq\emptyset}$ since in that case we can
use Lemma \ref{t-schauder} directly to obtain the desired result. We have
already shown that $u\in\mathcal{C}_{X}^{0,\alpha}(C_{R})$ so we are left with
the task to prove that%
\begin{equation}
||X_{i}u||_{\mathcal{C}_{X}^{0,\alpha}(C_{R})}=\sup_{C_{R}}|X_{i}%
u|+\sup_{C_{R}}\frac{|X_{i}u(x,t)-X_{i}u(\hat{x},\hat{t})|}{d_{p}%
((x,t),(\hat{x},\hat{t}))^{\alpha}}\leq\widehat{c}\label{sofia2}%
\end{equation}
and that%
\begin{equation}
\sup_{C_{R}}\left\vert X_{i}u(x,t)-X_{i}u(\tilde{x},\tilde{t})-\sum_{i=1}%
^{q}X_{i}u(\tilde{x},\tilde{t})(x_{i}-\tilde{x}_{i})\right\vert \leq
\widehat{c}d_{p}((x,t),(\tilde{x},\tilde{t}))^{1+\alpha}\label{sofia1}%
\end{equation}
for $i\in\{1,2,...,q\}$ and for some constant $\widehat{c}=\widehat
{c}(\mathcal{H},\alpha,f,\varphi)$.

\textbf{Step 1. }Notice that, as a consequence of Lemma \ref{lem1+} and the
fact that $u\geq$ $\varphi,$ we know that $X_{i}u(x,t)=X_{i}\varphi(x,t)$
whenever $(x,t)\in\mathcal{F}$. Now, let $(\hat{x},\hat{t})\in C_{2R}%
\cap\mathcal{F}$ and define the function
\begin{equation}
v(x,t):=u(x,t)-u(\hat{x},\hat{t})-\sum_{i=1}^{q}X_{i}u(\hat{x},\hat{t}%
)(x_{i}-\hat{x}_{i}).\label{nkelv}%
\end{equation}
We note that, in fact $v(x,t)=u(x,t)-P_{1}^{(\hat{x},\hat{t})}\varphi(x,t)$
and we may use Lemma \ref{lem1+} to deduce that%
\[
|v(x,t)|\leq c_{1}d_{p}((x,t),(\hat{x},\hat{t}))^{1+\alpha},
\]
for all $(x,t)\in C_{R}.$That is, (\ref{sofia1}) holds when the supremum is
taken over the set where $(x,t)\in C_{R}\backslash\mathcal{F}$, $(\hat{x}%
,\hat{t})\in C_{R}\mathcal{\cap F}$.

\textbf{Step 2.} Now we intend to prove that (\ref{sofia2}) also holds for
$(x,t)\in C_{R}\backslash\mathcal{F}$, $(\hat{x},\hat{t})\in C_{R}%
\mathcal{\cap F}$. To do this let $r:=d_{p}((x,t),\mathcal{F})$ and consider
the case when $d_{p}((x,t),(\hat{x},\hat{t}))\leq2r.$ Then, for $v$ defined as
in (\ref{nkelv}) we define the function%
\begin{equation}
w(y,s):=\frac{v^{r,(x,t)}(y,s)}{r^{1+\alpha}}\text{ \ on }C_{1}%
.\label{dubbelw}%
\end{equation}
This function satisfies the following on $C_{1/2}$, using the notation
introduced in Section \ref{SecFC},

\begin{description}
\item[$\mathit{i)}$] $|w|\leq c_{2}$ for some constant $c_{2},$ independent of
$r,$

\item[$\mathit{ii)}$] $\mathcal{H}^{r,(x,t)}w(y,s)=r^{1-\alpha}\left(
f^{r,(x,t)}-\sum_{i=1}^{q}b_{i}^{r,(x,t)}X_{i}\varphi(\hat{x},\hat{t})\right)
.$
\end{description}

\noindent By using the Schauder estimate in Lemma \ref{t-schauder} we see that%
\begin{equation}
||X_{i}w||_{\mathcal{C}_{X}^{0,\alpha}(C_{R})}\leq c_{2}^{\prime},\text{
\ \ }i\in\{1,2,...q\},\label{surf}%
\end{equation}
and since%
\begin{equation}
X_{i}w(0,0)=\frac{X_{i}u(x,t)-X_{i}u(\hat{x},\hat{t})}{r^{1+\alpha}%
},\label{surfa}%
\end{equation}
this implies the H\"{o}lder continuity of $X_{i}u$. Furthermore, we may
combine (\ref{surf}), (\ref{surfa}) and the fact that $X_{i}u(\hat{x},\hat
{t})=X_{i}\varphi(\hat{x},\hat{t})$ with $r=R$ to obtain%
\[
\sup_{C_{R}}|X_{i}u|\leq c_{2}^{\prime}+\sup_{C_{R}}|\varphi|.
\]
On the other hand, if $d_{p}((x,t),(\hat{x},\hat{t}))>2r$ and if $(\tilde
{x},\tilde{t})\in\mathcal{F}$ is such that $r=d_{p}((x,t),(\tilde{x},\tilde
{t})) $ we use the triangle inequality to obtain%
\begin{align*}
d_{p}((\hat{x},\hat{t}),(\tilde{x},\tilde{t}))  & \leq d_{p}((x,t),(\hat
{x},\hat{t}))+d_{p}((x,t),(\tilde{x},\tilde{t}))\\
& \leq3/2d_{p}((x,t),(\hat{x},\hat{t})).
\end{align*}
Hence,%
\begin{align*}
\left\vert X_{i}u(x,t)-X_{i}u(\hat{x},\hat{t})\right\vert  & \leq\left\vert
X_{i}u(x,t)-X_{i}u(\tilde{x},\tilde{t})\right\vert +\left\vert X_{i}%
u(\tilde{x},\tilde{t})-X_{i}u(\hat{x},\hat{t})\right\vert \\
& \leq d_{p}((x,t),(\tilde{x},\tilde{t}))^{\alpha}+\left\vert X_{i}%
\varphi(\tilde{x},\tilde{t})-X_{i}\varphi(\hat{x},\hat{t})\right\vert \\
& \leq c_{3}^{\prime}d_{p}((x,t),(\hat{x},\hat{t}))^{\alpha},
\end{align*}
where we also used the definition of $r$ and the fact that $X_{i}%
u(x,t)=X_{i}\varphi(x,t)$ whenever $(x,t)\in\mathcal{F}.$ That is, we have
proved that (\ref{sofia2}) holds when the supremum is taken over the set where
$(x,t)\in C_{R}\backslash\mathcal{F}$, $(\hat{x},\hat{t})\in\mathcal{F\cap
}C_{R}.$

\textbf{Step 3.} In this step, which completes the proof of (\ref{sofia2}), we
will prove that the statement also is valid for $(x,t),\ (\hat{x},\hat{t})\in
C_{R}\backslash\mathcal{F}.$ As before, we let $r$ and $(\tilde{x},\tilde
{t})\in$ $C_{2R}\cap\mathcal{F}$ be such that $r=d_{p}((x,t),\mathcal{F}%
)=d_{p}((x,t),(\tilde{x},\tilde{t}))$ and firstly we consider the case when
$(\hat{x},\hat{t})\in C_{r}(x,t)$. In analogy with (\ref{dubbelw}) we define a
function $w$ but we replace $(\hat{x},\hat{t})$ with $(x,t)$ in (\ref{nkelv}).
Using the same argument as above we find that (\ref{surf}) still holds and in
particular%
\[
|X_{i}w(y,s)-X_{i}w(0,0)|=r^{-\alpha}|X_{i}(u)((x,t)\circ\delta_{r}%
(y,s))-X_{i}(u)(x,t)|\leq c_{2}^{\prime}.
\]
If $(\hat{x},\hat{t})\in C_{R}\backslash C_{r}(x,t)$ we can use the triangle
inequality to obtain%
\[
|X_{i}u(x,t)-X_{i}u(\hat{x},\hat{t})|\leq|X_{i}u(x,t)-X_{i}u(\tilde{x}%
,\tilde{t})|+|X_{i}u(\tilde{x},\tilde{t})-X_{i}u(\hat{x},\hat{t})|.
\]
Now, since $(x,t),$ $(\hat{x},\hat{t})$ $\in C_{R}\backslash\mathcal{F}$ and
$(\tilde{x},\tilde{t})\in\mathcal{F\cap}C_{R}$ we can use the previous result
and%
\begin{align*}
|X_{i}u(x,t)-X_{i}u(\hat{x},\hat{t})|  & \leq\max\{c_{3},c_{3}^{\prime
}\}\left(  d_{p}((x,t),(\tilde{x},\tilde{t}))^{\alpha}+d_{p}((\tilde{x}%
,\tilde{t}),(\hat{x},\hat{t}))^{\alpha}\right) \\
& \leq c_{3}^{\prime\prime}d_{p}((x,t),(\hat{x},\hat{t}))^{\alpha},
\end{align*}
which concludes the proof of (\ref{sofia2}).

\textbf{Step 4.} It remains to prove that (\ref{sofia1}) holds for
$(x,t),\ (\hat{x},\hat{t})\in C_{R}\backslash\mathcal{F}$. Again, we let $r$
and $(\tilde{x},\tilde{t})\in\mathcal{F\cap}C_{R}$ be such that $r=d_{p}%
((x,t),\mathcal{F})=d_{p}((x,t),(\tilde{x},\tilde{t}))$ and we divide the
proof in two cases. Firstly we assume that $(\hat{x},\hat{t})\in
C_{R}\backslash C_{r/2}(x,t)$. Then we apply the triangle inequality and%
\begin{align*}
& \left\vert u(x,t)-u(\hat{x},\hat{t})-\sum_{i=1}^{q}X_{i}u(\hat{x},\hat
{t})(x_{i}-\hat{x}_{i})\right\vert \\
& \leq\left\vert u(x,t)-u(\tilde{x},\tilde{t})-\sum_{i=1}^{q}X_{i}u(\tilde
{x},\tilde{t})(x_{i}-\tilde{x}_{i})\right\vert +\left\vert u(\tilde{x}%
,\tilde{t})-u(\hat{x},\hat{t})-\sum_{i=1}^{q}X_{i}u(\hat{x},\hat{t})(\tilde
{x}_{i}-\hat{x}_{i})\right\vert \\
& +\sum_{i=1}^{q}|X_{i}(\tilde{x},\tilde{t})-X_{i}u(\hat{x},\hat{t}%
)|\cdot|\tilde{x}_{i}-x_{i}|\\
& \leq c_{4}^{\prime}\left(  d_{p}((x,t),(\hat{x},\hat{t}))^{1+\alpha}%
+d_{p}((\tilde{x},\tilde{t}),(\hat{x},\hat{t}))^{1+\alpha}+d_{p}((\tilde
{x},\tilde{t}),(\hat{x},\hat{t}))^{\alpha}d_{p}((\tilde{x},\tilde
{t}),(x,t))\right) \\
& \leq c_{4}^{\prime\prime}d_{p}((x,t),(\hat{x},\hat{t}))^{1+\alpha}.
\end{align*}
Now, assume that $(\hat{x},\hat{t})\in C_{r/2}(x,t)\cap C_{R}.$ We recall that
$u(x,t)-u(\hat{x},\hat{t})-\sum_{i=1}^{q}X_{i}u(\hat{x},\hat{t})(x_{i}-\hat
{x}_{i})$ is in fact $u(x,t)$ minus its Taylor polynomial of $\delta_{\lambda
}$-degree $1$ at $(\hat{x},\hat{t}).$ In this case we use that $u$ is a
solution to (\ref{e-obsobs}) and hence $u\in\mathcal{C}_{X}^{2,\alpha} $, see
Subsection \ref{secFL}. We can thus use Lemma \ref{t-schauder} to find
$\mathcal{C}_{X}^{2,\alpha}$-bounds on $u.$ Then we use these bounds in
Theorem \ref{Taylor} to obtain the desired result. Note that we also used that
we are working on a stratified group when using Theorem \ref{Taylor}. This
completes the proof in the case $m=1,$ $\gamma=1+\alpha$.

\subsection{Proof of Theorem \ref{th-int} when $m=2,$ $\gamma=2$}

To ease notation we assume, as in Lemma \ref{lemD}, that $\varphi\equiv0,$
note however that this is not restrictive. That is we assume that
\[
(u,f,g,\varphi)\in\mathcal{P}_{m}(\mathcal{H},\alpha,C,M_{1},M_{2},M_{3}),
\]
where $\varphi\equiv0$ and $M_{3}=0.$ As before we take $R\in(0,1)$ and
$C_{1}>1$ such that $C_{2R}\subseteq C_{2C_{1}R}(\hat{x},\hat{t})\subseteq C$
for all $(\hat{x},\hat{t})\in C_{2R}$ and we define $\mathcal{F}=\overline
{C}_{2R}\cap\{(x,t):u(x,t)=\varphi(x,t)\}.$ We intend to prove that there
exists a constant $\widehat{c}=\widehat{c}(\mathcal{H},\alpha,M_{1}%
,M_{2})<\infty$ such that%
\begin{equation}
||u||_{\mathcal{S}_{X}^{\infty}(C_{R})}\leq\widehat{c}.\label{pff}%
\end{equation}
For $(\hat{x},\hat{t})\in C_{R}\cap\left\{  (x,t):u(x,t)>0\right\}  $ we
define $\hat{r}=\hat{r}(\hat{x},\hat{t})=\sup\{r:C_{r}(\hat{x},\hat{t})\subset
C\cap\{(x,t):u(x,t)>0\}\}.$ By the maximum principle $\mathcal{F\cap\partial
}_{p}C_{\hat{r}}(\hat{x},\hat{t})\neq\emptyset$ which means that there exists
$(\tilde{x},\tilde{t})\in C_{2R}\cap\mathcal{F\cap\partial}_{p}C_{\hat{r}%
}(\hat{x},\hat{t})$ such that $C_{\hat{r}}(\hat{x},\hat{t})\subset
C_{\tilde{r}}(\tilde{x},\tilde{t})$ for some $\tilde{r},\ \hat{r}<\tilde
{r}<c_{1}\hat{r}$, for some universal constant $c_{1}.$ Now we use Lemma
\ref{lem1+} which states that for $r\in(0,\tilde{r})$ and $(x,t)\in
C_{r}(\tilde{x},\tilde{t})\cap C,$%
\begin{equation}
|u(x,t)|\leq c_{2}r^{2}.\label{uhh}%
\end{equation}
Thereafter we define, for $(x,t)\in C$,%
\[
v(x,t):=\frac{u^{\hat{r},(\hat{x},\hat{t})}(x,t)}{\hat{r}^{2}},
\]
and in particular there holds%
\[
\mathcal{H}^{\hat{r},(\hat{x},\hat{t})}v=f^{\hat{r},(\hat{x},\hat{t})}\text{
\ \ in }C.
\]
For notation, see Subsection \ref{SecFC}. Furthermore, by (\ref{uhh}) and
(\ref{fredag}),%
\[
||v||_{L^{\infty}(C)}\leq c_{2}\text{ \ \ and \ \ \ }||f^{\hat{r},(\hat
{x},\hat{t})}||_{\mathcal{C}_{X}^{0,\alpha}(C)}\leq M_{2}.
\]
In particular, we can use Lemma \ref{t-schauder} to conclude that
\[
||v||_{\mathcal{S}_{X}^{\infty}(C_{1/2}(\hat{x},\hat{t}))}\leq
||v||_{\mathcal{C}_{X}^{2,\alpha}(C_{1/2}(\hat{x},\hat{t}))}\leq c,
\]
and (\ref{pff}) follows immediately. This completes the proof of Theorem
\ref{th-int}.\hfill$\Box$

\section{Proof of Theorem \ref{th}}

The proof is very similar to the proof of Theorem \ref{th-int} but differs in
that we, near the initial state, have to rely on Lemma \ref{lem} instead of
Lemma \ref{lem1+}. We will give a detailed proof of part $\mathit{i)}$ of
Theorem \ref{th}, leaving part $\mathit{ii)}$ and $\mathit{iii)}$. For this
reason, let $R\in(0,1)$ and $C_{1}>1$ be such that $C_{2R}\subseteq
C_{2C_{1}R}(\hat{x},\hat{t})\subseteq C$ for all $(\hat{x},\hat{t})\in
C_{2R}.$ We will now study the regularity of $u$ near the initial state and we
note that, by Subsection \ref{SecFC}, it is enough to consider $\Omega=C^{+}$
and $\Omega^{\prime}=C_{R}^{+}$. To prove part $\mathit{i)}$\textit{\ }we need
to show that there exists a constant $c=c\left(  \mathcal{H},\alpha
,||f||_{\mathcal{C}_{X}^{0,\alpha}(C^{+})},||g||_{L_{{}}^{\infty}(C^{+}%
)},||\varphi||_{\mathcal{C}_{X}^{0,\alpha}(C^{+})}\right)  $ such that%
\begin{equation}
\underset{(x,t)\neq(\hat{x},\hat{t})}{\sup_{(x,t),(\hat{x},\hat{t})\in
C_{R}^{+}}}\frac{|u(x,t)-u(\hat{x},\hat{t})|}{d_{p}((x,t),(\hat{x},\hat
{t}))^{\alpha}}\leq c.\label{betti1}%
\end{equation}
If $\hat{t}=0$ in (\ref{betti1}) we use Lemma \ref{lem} on $C_{d_{p}%
((x,t),(\hat{x},0))}^{+}(\hat{x},0)$ with $m=0$ and obtain%
\begin{align}
|u(x,t)-u(\hat{x},0)|  & \leq|u(x,t)-g(x,t)|+|g(x,t)-g(\hat{x},0)|\nonumber\\
& \leq c_{1}d_{p}((x,t),(\hat{x},0))^{\alpha},\label{betti2}%
\end{align}
since $g\in C(\overline{C^{+}})$. Hence (\ref{betti1}) is valid if either one
of $t,\hat{t}$ vanishes and hereafter we assume that both $t$ and $\hat{t}$
are strictly positive. In that case we note that $d_{p}((x,t),(x,0))=\sqrt{t}$
and we will divide the rest of the proof in two cases.

\textbf{Case 1: }We assume that $(\hat{x},\hat{t})\in C_{R}^{+}\backslash
C_{\sqrt{t}/2}(x,t)$ which implies that%
\begin{align*}
d_{p}((\hat{x},\hat{t}),(x,0))  & \leq d_{p}((\hat{x},\hat{t}),(x,t))+d_{p}%
((x,t),(x,0))\\
& \leq3d_{p}((\hat{x},\hat{t}),(x,t))
\end{align*}
and%
\[
d_{p}((x,t),(x,0))\leq2d_{p}((\hat{x},\hat{t}),(x,t)).
\]
Therefore, by (\ref{betti2}),%
\begin{align*}
|u(x,t)-u(\hat{x},\hat{t})|  & \leq|u(x,t)-u(x,0)|+|u(x,0)-u(\hat{x},\hat
{t})|\\
& \leq d_{p}((x,t),(x,0))^{\alpha}+d_{p}((\hat{x},\hat{t}),(x,0))^{\alpha}\\
& \leq5d_{p}((\hat{x},\hat{t}),(x,t))^{\alpha},
\end{align*}
which closes case 1.

\textbf{Case 2:} We assume that $(\hat{x},\hat{t})\in C_{R}^{+}\cap
C_{\sqrt{t}/2}(x,t)$ and note that $C_{\sqrt{t}}(x,t)\subset C_{2R}^{+}.$ By
(\ref{betti2})%
\begin{align}
||u-u(x,t)||_{L^{\infty}\left(  C_{\sqrt{t}}(x,t)\right)  }  & \leq
||u-u(x,0)||_{L^{\infty}\left(  C_{\sqrt{t}}(x,t)\right)  }%
+|u(x,t)-u(x,0)|\nonumber\\
& \leq c_{2}\left(  \sqrt{t}\right)  ^{\alpha}\label{betti3}%
\end{align}
and we define
\[
v(y,s):=\frac{u^{\sqrt{t},(x,t)}(y,s)-u(x,t)}{\left(  \sqrt{t}\right)
^{\alpha}}\text{ \ \ for }(y,s)\in C\text{.}%
\]
Using (\ref{betti3}) we see that $||v||_{L^{\infty}(C)}\leq c_{2}$ and by
using properties of $v$ (see Subsection \ref{SecFC}) and Theorem \ref{th-int}
we see that%
\[
|v(y,s)|\leq c_{3}d_{p}((0,0),(y,s))^{\alpha}\text{ \ \ for }(y,s)\in C_{1/2}.
\]
By definition of $v$ this is equivalent to
\[
|u(x,t)-u(\hat{x},\hat{t})|\leq cd_{p}((x,t),(\hat{x},\hat{t}))^{\alpha},
\]
and this completes the proof of part $\mathit{i)}$ of Theorem \ref{th}. As
previously mentioned the rest of the proof, i.e., the proof of $\mathit{i)}$
and $\mathit{ii),}$ is much in line with the corresponding statements in
Theorem \ref{th-int} and we omit the details. This completes the proof.
\hfill$\Box$

\section{Generalizations, further developments and open problems}

When we proved the existence of strong solutions to the obstacle problem
(\ref{e-obs}) in \cite{FGN} we were able to carry out the proofs using less
restrictive assumptions than in the present paper. The main difference is
that, unlike in \cite{FGN}, we here assume that the vector fields
$\{X_{1},...,X_{q}\}$ are generators of the first layer of a stratified,
homogeneous group. In addition we assume that $\{X_{1},...,X_{q}\}$ are left
invariant and homogeneous of degree one. In \cite{FGN} we did not have to
restrict ourselves to this case and the main tool for carrying out the proofs
was the lifting-approximation technique of Rotschild and Stein \cite{RS}. The
main difficulties to overcome in the present paper, if one should try to relax
these assumptions, are that of polynomial approximations and scaling, see
Section \ref{SecFS} and Section \ref{SecFC} respectively. It may be possible
to use the lifting-approximation technique here as well, however, it is not
clear if $\mathcal{C}_{X}^{m,\alpha}$ estimates carry through the
lifting-approximation machinery. A partial (affirmative) answer is given in
\cite{BB07} which states that if a function $u$ is lifted to $\widetilde{u}$
and if $\widetilde{u}\in\mathcal{C}_{\widetilde{X}}^{m,\alpha}$, where
$\widetilde{X}$ is the lifted vector field, then $u\in\mathcal{C}_{X\text{ }%
}^{m,\alpha}$ as well. Whether or not a similar estimate holds for the
approximation part is unclear.

Another generalization would be to consider the operator
\begin{equation}
\mathcal{H}=\sum_{i,j=1}^{q}a_{ij}(x,t)X_{i}X_{j}+\sum_{i=1}^{q}%
b_{i}(x,t)X_{i}+X_{0},\ \ (x,t)\in\mathbf{R}^{n+1}.\label{operator1}%
\end{equation}
When we began our study in \cite{FGN} certain estimates for the operator in
(\ref{operator1}) were missing, in particular Schauder type estimates,
interior $\mathcal{S}^{p}$-estimates and embedding theorems corresponding to
Theorems 1.2-1.4 in \cite{FGN}. However, for operators (\ref{operator}) at
least the necessary Schauder type estimates were established and so we chose
to work with (\ref{operator}) instead of (\ref{operator1}). A recent preprint
by Bramanti and Zhu \cite{BZ11} has made two of the crucial estimates
available also for operators (\ref{operator1}). Two prove the embedding
theorem we need estimates on the fundamental solution. If the vector fields
$X_{0},X_{1},...,X_{q}$ are left invariant on a homogeneous group and if
$X_{1},...,X_{q}$ are homogeneous of degree $1$ while $X_{0}$ is homogeneous
of degree $2$ then we may use results of Folland \cite{F75}, in the more
general case this is an open question. Provided with this estimate, and after
proving the existence of strong solutions to the obstacle problem, the same
method we have presented here will most likely be usable for the operator
(\ref{operator1}).

Finally, it would be interesting to study the regularity of the free boundary
in the setting studied in the present paper.

\end{document}